\documentclass[lettersize,journal]{IEEEtran}
\usepackage{amsmath,amssymb,amsfonts}
\usepackage[ruled,vlined]{algorithm2e}
\usepackage{graphicx}
\usepackage{textcomp}
\usepackage{latexsym}
\usepackage{graphicx,subfigure}
\usepackage[hidelinks]{hyperref}
\usepackage{fancyhdr}
\usepackage{bm}
\usepackage{cite}
\usepackage{mathrsfs}
\usepackage{subeqnarray}
\usepackage{cases}
\usepackage{arydshln}
\usepackage{algpseudocode}
\usepackage{multirow}

\allowdisplaybreaks[4]

\newtheorem{theorem}{Theorem}
\newtheorem{lem}{Lemma}
\newtheorem{remark}{Remark}
\newtheorem{ass}{Assumption}
\newtheorem{definition}{Definition}
\newtheorem{prop}{Proposition}

\begin{document}
	
\title{\Huge Extended Zero-Gradient-Sum Approach for Constrained Distributed Optimization with Free Initialization}
% With User-Defined Convergence Performance }

\author{Xinli~Shi,~\IEEEmembership{Senior Member, IEEE},  Xinghuo~Yu,~\IEEEmembership{Fellow,~IEEE}, Guanghui~Wen,~\IEEEmembership{Senior Member, IEEE}, and Xiangping Xu,~\IEEEmembership{Member, IEEE}
	%\author{
%\thanks{This work was supported by the National Natural Science Foundation of China under Grant Nos. 62203108 and 62073079, and the Australian Research Council under Discovery Program DP200101199, the General Joint Fund of the Equipment Advance Research Program of Ministry of Education under Grant No. 8091B022114, Jiangsu Province Excellent Postdoctoral Program under Grant No. 2022ZB131, and China Postdoctoral Science Foundation under Grant Nos. 2022M720720 and 2023T160105.}
	\thanks{X.~Shi is with the School of Cyber Science and Engineering, Southeast University, Nanjing 210096, China (e-mail: xinli\_shi@seu.edu.cn).}
	\thanks{X. Yu is with the School of Engineering, RMIT University, Melbourne, VIC 3001, Australia (e-mail: x.yu@rmit.edu.au).}
	\thanks{G. Wen and X. Xu are with the School of Mathematics, Southeast University, Nanjing 210096, China (e-mail: ghwen@seu.edu.cn; xpxu2021@seu.edu.cn).}
}

%	\thanks{G. Wen and J. Cao are with the School of Mathematics, Southeast University, Nanjing 210096, China (e-mail: wenguanghui@gmail.com;jdcao@seu.edu.cn).}
	%	\thanks{J. Cao is with the School of Mathematics, Frontiers Science Center for Mobile Information Communication and Security, Southeast University, Nanjing 210096, China; Yonsei Frontier Lab, Yonsei University, Seoul 03722, South Korea (e-mail: jdcao@seu.edu.cn).}

\maketitle
\thispagestyle{fancy}
\fancyhead{}
\chead{\large This work has been submitted to the IEEE for possible publication. Copyright may be transferred without notice, after which this version may no longer be accessible.}

\begin{abstract}
This paper proposes an extended zero-gradient-sum (EZGS) approach for solving constrained distributed optimization (DO) with free initialization. A Newton-based continuous-time algorithm (CTA) is first designed for general constrained optimization and then extended to solve constrained DO based on the EZGS method. It is shown that for typical consensus protocols, the EZGS CTA can achieve the performance with exponential/finite/fixed/prescribed-time convergence. Particularly, the nonlinear consensus protocols for finite-time EZGS algorithms can have heterogeneous power coefficients. The prescribed-time EZGS dynamics is continuous and uniformly bounded, which can achieve the optimal solution in one stage. Moreover, the barrier method is employed to tackle the inequality constraints. Finally, the performance of the proposed algorithms is verified by numerical examples.
\end{abstract}
\begin{IEEEkeywords}
Zero-gradient-sum method, constrained distributed optimization, convergence rate, finite-time stability, barrier method 
\end{IEEEkeywords}

%% main text

\section{Introduction}
The past decade has witnessed rapid development of distributed optimization (DO) in both theory and applications. From the time scale, the existing DO algorithms can be categorized into discrete- and continuous-timeh as robotics and unmanned vehicles with continuous-time physical dynamics \cite{yang2019survey}. Moreover, the theoretical a types. With the development of cyber-physical systems, the continuous-time algorithm (CTA) is preferred in some practical systems sucnalysis of CTAs is more convenient based on Lyapunov stability methods. By transformation or discretization, CTA can be converted to the discrete-time algorithm (DTA). Compared with DTAs, some CTAs can achieve fascinating performance, such as finite-time (FT) convergence based on FT stability theory \cite{Bhat2000}. From the types of problems, DO can be classified into unconstrained and constrained optimization. The existing CTAs for unconstrained/constrained DO includes primal-dual methods \cite{gharesifard2013TAC,Guo2022TNNLS, Zhu2020TC}, neurodynamic approaches \cite{Jiang2021,Xia2023,Zhao2023}, zero-gradient-sum (ZGS) algorithms \cite{Lu2012TAC, chen2016event,song2016finite, Wu2021},  just to name a few. For more recent review and applications, one can refer to the survey \cite{yang2019survey} and references therein. 

%\cite{Lin2016TAC, song2016finite, Shi2022TAC,Shi2023Tcyber}

In the aforementioned works, most algorithms can only achieve asymptotic or exponential (EXP) convergence. To achieve FT convergence, some variants of existing algorithms have been proposed \cite{song2016finite, Wu2021, Shi2023CSL,Shi2022TAC,Shi2023Tcyber}. 
%In \cite{Lin2016TAC}, a discontinuous FT algorithm is provided for DO with cost functions subject to quadratic growth. 
By combining the ZGS scheme \cite{Lu2012TAC} and nonlinear consensus protocols, FT ZGS dynamics are provided in \cite{song2016finite,Wu2021} for solving unconstrained DO. The advantage of ZGS methods is that only the primal states are involved; however, to ensure the ZGS property, the initial local states are required to be the minimizers of local cost functions. To achieve free initialization, an extended ZGS dynamics is provided for solving time-varying DO \cite{Shi2023CSL}. Based on FT distributed average tracking (DAT) protocols and FT gradient dynamics (FT-GD), several discontinuous CTAs with globally bounded inputs are proposed in \cite{Shi2022TAC} for both unconstrained and constrained DO. The convergence analysis can be summarized in three-stages: FT-consensus, FT-DAT estimation and FT-GD. For the optimization with explicit equation and inequality constraints, based on nonsmooth analysis, an FT convergent primal-dual gradient dynamics (FT-PDGD) is proposed in \cite{Shi2023Tcyber} and further applied to solve DO with explicit constraints. 

To achieve finite settling time independent of the initial states, several fixed-time (FxT) algorithms have been provided for various DO problems \cite{Wang2020auto,Ma2023TC,Garg2020cdc}. For example, a two-stage FxT algorithm is considered in \cite{Wang2020auto} and \cite{Ma2023TC} for unconstrained DO based on the FxT-GD and ZGS method; however, the settling time estimation in the second stage is related to parameters determined by the initial states as pointed out in \cite{Wu2021}. A three-stage FxT algorithm is presented in \cite{Garg2020cdc} for unconstrained DO based on FxT-DAT estimation, FxT-consensus and FxT-GD. FxT consensus is considered in the DO dynamics in \cite{Ning2019TC}. In \cite{Guo2023}, a sliding mode-based ZGS algorithm is provided for unconstrained DO with FxT convergence and free initialization. Actually, the methods employed for achieving the free local initialization in \cite{Guo2023} and \cite{Shi2023CSL} are equivalent. Moreover, when the settling time can be pre-assigned by a designer, which is independent of the initial states and topology, the prescribed-time (PT) convergence is reached. For example, similar to \cite{Garg2020cdc}, the PT DO is investigated with three cascading stages in \cite{gong2021dis} based on PT-DAT estimation, PT-consensus and PT-GD, where the global average gradient and Hessian need to be estimated by DAT dynamics first with high communication/computation cost. Multiobjective DO 
with PT approximate convergence is considered in \cite{Liu2023}.

Although several algorithms have been proposed based on ZGS method, there still exist three limitations: 1) all the existing ZGS-based algorithms are for solving unconstrained DO; 2) to ensure the ZGS property, the local minimizers are required initially or achieved by another FxT-GD; 3) for FT/FxT ZGS methods \cite{song2016finite,Wu2021,Wang2020auto}, the nonlinear consensus protocols are not fully distributed since the power coefficients are homogeneous. For DAT based methods (e.g., \cite{Shi2022TAC,Garg2020cdc,gong2021dis}), as the estimation of the global average gradient/Hessian is required, it causes higher communication burden compared with ZGS based methods. For FT-PDGD \cite{Shi2023Tcyber}, an additional condition is imposed on the inequality constraints and it can only achieve FT convergence. Besides, the right-hand sides of the CTAs in \cite{gong2021dis,Shi2022TAC,Shi2023Tcyber} are all discontinuous, increasing the difficulty in the distributed implementation. 

%Motivated by the above issues, in this work, we provide an extended ZGS (EZGS) framework for solving constrained DO with user-defined performance. That is to say, with typical consensus protocols, the optimal primal-dual solution can be achieved with exponential/FT/FxT/PT convergence. 
%Moreover, the proposed FT/FxT EZGS algorithms are fully distributed with heterogeneous coefficients. 
%To eliminate the special initialization, we introduce an auxiliary dynamics to ensure the ZGS property and then the local states converge to the optimal solution. To tackle the inequality constraints, the barrier method is employed using log-barrier penalty functions motivated by \cite{fazlyab2017pred}. 
%Besides, the right-hand sides of the CTAs in \cite{Shi2020TAC,Lin2016TAC,gong2021dis,Shi2023Tcyber} are all discontinuous, increasing the difficulty in the distributed implementation. 

%To summarize, among the existing FT/FxT/PT convergent algorithms are mainly designed for unconstrained DO with special initialization \cite{song2016finite,Wu2021}, or quadratic-like objective functions \cite{Lin2016TAC,Feng2020}, or using DAT estimation with several stages and high communication cost \cite{Lin2016TAC,Shi2022TAC,Garg2020cdc,gong2021dis}. 
To summarize, there are few works on the FT/FxT/PT algorithms for general constrained DO with free initialization. Motivated by the above issues, this work is to fill the gap and the contributions can be summarized as follows.  

(1) We first develop a Newton-based CTA for solving equation constrained optimization in a unified framework, which can achieve the EXP/FT/FxT/PT convergence with typical protocols, extending \cite{Garg2021tac} (resp. \cite{Shi2023Tcyber}) where only FxT (resp. FT) convergence can be obtained by applying FxT gradient flows (resp. FT-PDGD). Using log-barrier penalty functions for inequality constraints, the primal constrained optimization is transformed into an approximate one. By the modified CTA with a time-varying penalty parameter, the state can track the time-varying approximate solution in FT/FxT/PT and then converges to the optimal solution to the primal problem exponentially. 

(2) Based on the proposed Newton-based CTA, we propose an extended ZGS (EZGS) approach for designing distributed CTAs to solve constrained DO with general constraints and free initialization, compared with the existing ZGS methods \cite{Lu2012TAC, chen2016event,song2016finite, Wu2021,Guo2023,Shi2023CSL} only targeted for unconstrained DO. Motivated by \cite{Shi2023CSL}, an auxiliary dynamics is introduced in EZGS approach to ensure the final ZGS property from any initialization. The EZGS approach can achieve the EXP/FT/FxT/PT convergence with typical consensus protocols. Particularly, the nonlinear consensus protocols for FT/FxT EZGS algorithms have heterogeneous power coefficients, compared with \cite{song2016finite,Wu2021,Wang2020auto,Guo2023}. The dynamics of the proposed PT EZGS algorithm is shown to be continuous and uniformly bounded, and the PT convergence can be achieved in one stage, compared with \cite{Wang2020auto,Garg2020cdc,gong2021dis,Ma2023TC}. The proposed EZGS approach is further extended to solve DO with inequalities using log-barrier penalty functions, which can achieve EXP/FT/FxT/PT convergence to an approximate optimal solution. 

%(3) We further give a discrete-time version of ZGS method for solving constrained DO over UJSC and balanced digraphs. The stepsize condition is derived to ensure the convergence. Motivated by the central path-following method (see \cite[Alg. 11.1]{boyd2004convex}), we propose a DPF-ZFS algorithm for solving DO with general inequality constraints, which can be regarded as a distributed version of the interior-point method.

The remainder of this work is organized as follows. Section \ref{section-2} gives the preliminary notations used in the later parts. The constrained DO problem statement and centralized algorithm are provided in Section \ref{Problem}. The proposed distributed EZGS dynamics and convergence analysis are shown in Section \ref{Sec-EZGS}. Finally, the numerical example is conducted in Section \ref{Numberical}. Conclusions are drawn in Section \ref{conclusion}.

%\vspace{-0.3cm}
\section{Preliminaries}\label{section-2}
\subsection{Notation and Network Representation}
Let $\mathbb{R}_+^n$ be the set of $n$-dimensional non-negative vectors. We use $\mathbf{1}_n$ $\in \mathbb{R}^n$ and $I_n\in \mathbb{R}^{n\times n}$ to represent the vector of all ones and identity matrix, respectively. Define $\langle n\rangle=\{1,2,\cdots,n\}$. For $x=[x_{1},\ldots,x_{n}]^{T}\in \mathbb{R}^n$, the Euclidean and 1-norm of $x$ are denoted by $\|x\|$ and $\|x\|_1$, respectively; $\text{sgn}^{\alpha}(x)\triangleq[\text{sign}(x_1)|x_1|^{\alpha},...,\text{sign}(x_n)|x_n|^{\alpha}]^T$ for $\alpha\in  \mathbb{R}_+$, and $\text{sgn}^{[\alpha]}(x)\triangleq[\text{sign}(x_1)|x_1|^{\alpha_1},...,\text{sign}(x_n)|x_n|^{\alpha_n}]^T$ for $\alpha =[\alpha_i]_{i\in \langle n\rangle}\in \mathbb{R}_+^n$, where $\text{sign}(\cdot)$ is the signum function. By default, we write $\text{sgn}(x)=\text{sgn}^{0}(x)=\text{sign}(x)$. For a matrix $M\in \mathbb{R}^{n\times n}$, $\mathcal{I}(M)$ and $\mathcal{N}(M)$ denotes the image and null space of the linear map $M$, respectively. When $M$ is positive semidefinite, $\lambda_2(M)$ is the smallest positive eigenvalue of $M$. 
A continuously differentiable function $f(x): \mathbb{R}^n \rightarrow \mathbb{R}$ is called
$\theta$-strongly convex if there exists a constant $\theta>0$ such that for any $x,y \in \mathbb{R}^n$, 
$f(y)-f(x) -\nabla f(x)^T(y-x) \geq \frac{\theta}{2}\|y-x\|^2$.

An undirected graph is represented by $\mathcal{G}(\mathcal{V},\mathcal{E},A_0)$, where $\mathcal{V}=\langle N\rangle$ and $\mathcal{E}\subseteq \mathcal{V} \times \mathcal{V}$ denote the node and edge set, respectively. The adjacency $A_0=[a_{ij}]_{N\times N}$ is symmetric and defined as: $a_{ij}=1$ iff $(j,i)\in \mathcal{E}$ and $a_{ii}=0, \forall i\in \mathcal{V}$. The \textit{Laplacian matrix} is given by $L(\mathcal{G})=\text{diag}(A_0\bm{1}_N)-A_0$. The incidence matrix of $\mathcal{G}$ is represented by $B=[b_{ik}]_{N \times m}$ with $b_{ik}=-b_{jk}= a_{ij}$ if $i<j$ for any $e_k=\{i,j\} \in \mathcal{E}$, where $m$ is the number of undirected edges of $\mathcal{G}$. 

%A time-varying graph $\mathcal{G}(k)$ is said to be UJSC if there exists $T$ > 0 such that for any tk, the union ∪t∈[tk,tk+T]G(t) is strongly connected

%Moreover, one can define a digraph $\mathcal{G}_d(\mathcal{V}, \mathcal{E}_d)$ corresponding to $\mathcal{G}(\mathcal{V},\mathcal{E},A_0)$ by assigning an orientation from $j$ to $i$ if $i<j$ for each $e_k=(i,j)\in \mathcal{E}$, and define the weighted incidence matrix $B=[b_{ik}]_{N \times m}$ correspondingly with
%\begin{align}\label{incidence}
%	b_{ik} = \left\{
%	\begin{array}{ll}
%	     a_{ij}, & \text{if} \  i<j,\\
%		-a_{ij}, & \text{if} \  i>j. \\
%	\end{array}
%	\right.
%\end{align}

\subsection{EXP/FT/FxT/PT Stability}
Consider the following differential system
	\begin{align}\label{auto_dynamic}
		\dot{x}(t)=f(x(t), t), x(0)=x_0
	\end{align}
with $x\in \mathbb{R}^n$. The concept of the EXP/FT/FxT/PT stability are provided in Definition \ref{def-FTC}. The related criterion can refer to \cite[Lem. 1]{Shi2020auto} and Lemma \ref{lem-PT}. When $f$ is discontinuous, the Filippov solution of \eqref{auto_dynamic} will be used, which is an absolutely continuous map $X:I\subset \mathbb{R} \rightarrow \mathbb{R}^n $ satisfying the differential inclusion
\begin{align}\label{filippov}
	\dot{X}(t) \overset{a.e.}\in \mathcal{F}[f](X(t),t), \ \forall t \in I.
\end{align}
The Filippov set-valued map $\mathcal{F}[f]: \mathbb{R}^n \times \mathbb{R}_+ \rightarrow 2^{\mathbb{R}^n}$ is defined by 
\begin{align}\label{setmap}
	\mathcal{F}[f](x,t) \triangleq  \bigcap_{\delta>0} \bigcap_{\omega(\mathcal{N}_0)=0}\overline{co}\{f(\hat{x},t) : \hat{x} \in \mathcal{B}_{\delta}(x)\backslash \mathcal{N}_0 \}, 
\end{align} 
where $\mathcal{B}_{\delta}(x)$ is an open ball centered at $x$ with radius $\delta>0$, and $\omega(\mathcal{N}_0)$ denotes the Lebesgue measure of $\mathcal{N}_0$.  
\begin{definition}
A sequence $\{x(t)\}_{t\geq0 }$ generated by \eqref{auto_dynamic} is called FT convergent to a point $\bar{x}$ if there exists a finite time $T(x_0)$ such that 
\begin{align*}%\label{finite-time}
		\lim_{t \rightarrow T(x_0)}x(t; x_0) =\bar{x} \ \text{and} \ x(t)=\bar{x}, \forall t\geq T(x_0).
\end{align*}
If the above $T(x_0)$ is uniformly bounded, i.e., $\exists T_{\max}<+\infty, T(x_0)\leq T_{\max}, \forall x_0 \in \mathbb{R}^n$, it is called FxT convergent. If such $T_{\max}$ can be preassigned, we call it PT convergent. 
\end{definition}

\begin{definition}\label{def-FTC}
Let $x(t; x_0)$ be any solution starting from $x_0 \in \mathbb{R}^n$. The origin of the system (\ref{auto_dynamic}) is \textit{exponentially stable} with the decay rate $r_0>0$ if it satisfies that $\|x(t; x_0)\|\leq \gamma(\|x_0\|)e^{-r_0 t}, \forall t\geq 0$, where $\gamma: \mathbb{R}_+\rightarrow  \mathbb{R}_+$. The origin is called globally \textit{FT/FxT/PT stable} if it is Lyapunov stable and $x(t)$ is FT/FxT/PT convergent to 0. 
\end{definition}
%In this paper, the finite/fixed-time stability of \eqref{auto_dynamic} given in Definition \eqref{def-FTC} will be investigated, for which the following lemmas are helpful.
%\begin{lem}\label{lem-FTC}
%	Let $x(t; x_0)$ be any solution of the system (\ref{auto_dynamic}) starting from $x_0\in  \mathbb{R}^n$. Suppose that there exists a continuous radially unbounded function $V: \mathbb{R}^n\rightarrow \mathbb{R}_+$ satisfies $V(x)=0 \Rightarrow x=0$. 
%	\begin{enumerate}
%		\item \cite[Thm. 4.2]{Bhat2000} If $\dot{V}(x(t))\leq -cV^r(x(t))$ with $c>0, r\in (0,1)$, the origin is FT stable with settling time $T(x_0)\leq \frac{V^{1-r}(x(0))}{c(1-r)}$.
%		\item \cite[Lem. 1]{Polyakov2012} If $V(x)$ satisfies that
%		\begin{align}
%			\dot{V}(x(t))\leq -(c_1 V^{r_1}(x(t))+ c_1 V^{r_2}(x(t) )^k
%		\end{align}
%		for some positive scalars $c_1,c_2,r_1,r_2: 0<r_1k<1, r_2k>1$, the origin is FxT stable with settling time estimated by
%		\begin{align}
%			T(x_0)\leq \frac{1}{c_1^k(1-r_1k)} + \frac{1}{c_2^k(r_2k-1)}.
%		\end{align}
%	\end{enumerate}
%\end{lem}

%\begin{lem}\label{lem-ineq}
%	For $z_1,z_2,\cdots,z_N \geq 0$, it holds that
%	\begin{align}
%		&\sum_{i=1}^{N }z_i^l \geq (\sum_{i=1}^{N }z_i)^l, \quad \text{if} \quad l \in (0,1) \\
%		&\sum_{i=1}^{N }z_i^l \geq N^{1-l}(\sum_{i=1}^{N }z_i)^l, \quad \text{if} \quad l >1. 
%	\end{align}
%\end{lem}

\section{Problem statement and Approach}\label{Problem}
\subsection{Problem Statement}\label{sec2-A}
Consider the following equality constrained DO with $N$ agents over an undirected network $\mathcal{G}(\mathcal{V}, \mathcal{E})$
\begin{align}\label{CO}
	\min_{x\in \mathbb{R}^n} \sum_{i=1}^Nf_i(x),\ s.t. \ A_ix =b_i, \ g_i(x)\leq 0, \ \forall i\in \mathcal{V}
\end{align}
%\begin{align}
%	\min_{x\in \mathbb{R}^n} \sum_{i=1}^Nf_i(x), \ s.t. & \ A_ix =b_i, g_i(x)\leq 0, \ i=1,\cdots, N 
%\end{align}
where $f_i: \mathbb{R}^n  \rightarrow \mathbb{R}$ is the local cost function, $A_ix =b_i$ with $A_i \in \mathbb{R}^{r_i\times n}$ and $g_i(x)\leq 0$ are the equality and inequality constraints, only known by the agent $i$. Suppose that $g_i: \mathbb{R}^n  \rightarrow \mathbb{R}^{p_i}$ is convex and twice continuously differentiable. %When there exists no equation constraint at agent $i$, we suppose that $A_i=0$. 
%In the latter parts, we denote $\mathcal{V}_1 = \{i\in \mathcal{V}: A_i\neq 0\}$ and $\mathcal{V}_2 = \mathcal{V}\backslash \mathcal{V}_1$, representing the set of agents with and without equation constraints.   
Moreover, the optimization \eqref{CO} is supposed to satisfy Assumptions \ref{ass-f}-\ref{ass-g}, which indicates that there exists a unique solution to \eqref{CO} denoted by $x^*$. By Assumption \ref{ass-g}, the optimal primal-dual solution can be determined by Karush–Kuhn–Tucker (KKT)
conditions \cite[Sec. 5]{boyd2004convex}. When the inequalities are absent, we denote the set of optimal primal-dual points to \eqref{CO} as $\mathcal{Z}^*=\{x^*\}\times \mathcal{D}^*$ with $\mathcal{D}^*$ representing the set of optimal dual solutions corresponding to the equality. This work aims to design distributed algorithms with only local information that solve \eqref{CO} with expected convergence performance, such as EXP/FT/FxT/PT convergence. 

\begin{ass}\label{ass-f} (Strong Convexity Condition)
	Each local function $f_{i}(x), i\in \mathcal{V} $ is twice continuously differentiable and $\theta_i$-strongly convex. 
\end{ass}

\begin{ass}\label{ass-A} (Full Rank Condition)
Each matrix $A_i, i\in \mathcal{V} $ has full row rank, and $\cap_{i\in \mathcal{V}} \mathcal{S}_i\neq \emptyset$ with $\mathcal{S}_i =\{x_i|A_ix_i=b_i\}$. 
\end{ass}

\begin{ass}\label{ass-g} (Slater’s Condition)
There exists at least one interior point $x'$ of the feasible set, i.e., 
$g_i(x')<0, A_ix'=b_i, \forall i\in \mathcal{V}$. 
\end{ass}

\subsection{Newton-based CTA}
We first consider the following optimization
\begin{align}\label{CO-cen}
	\min_{x\in \mathbb{R}^n} F(x), \ s.t.  \ Ax =b, g(x)\leq 0,
\end{align}
where $F(x) \triangleq \sum_{i=1}^Nf_i(x)$ and $A=[A_i]_{i\in \mathcal{V}} \in \mathbb{R}^{r\times n}, b =[b_i]_{i\in \mathcal{V}} \in \mathbb{R}^r$ and $g(x)=[g_i(x)]_{i\in \mathcal{V}}: \mathbb{R}^n  \rightarrow \mathbb{R}^{p}$ with $r=\sum_{i=1}^{N}r_i$ and $p=\sum_{i=1}^{N}p_i$. Let $\mathcal{X} = \{x \in \mathbb{R}^n| Ax =b, g(x)\leq 0\}$ denote the feasible state set. In the following, We aim to provide a continuous-time dynamics that solves \eqref{CO-cen} with expected convergence performance.

First, we consider the optimization \eqref{CO-cen} without inequality constraints and introduce the Lagrangian function as $\mathcal{L}(x,\lambda)=F(x)+\lambda^T(Ax-b) $, where $\lambda$ is the dual variable related to the equation constraint. Let $z=(x,\lambda)$. When $F(x)$ is strongly convex and $A$ is of full row-rank, $\nabla^{2} \mathcal{L}(z)$ is invertible and there exists a unique optimal primal-dual point $z^*=(x^*,\lambda^*)$. Under this condition, the proposed dynamics is given by
\begin{subequations}\label{prot-cent}
	%\small 
	\begin{numcases}{}
		\dot{z} =-(\nabla^{2} \mathcal{L}(z))^{-1} \varphi(y,t), \ \ \quad \label{prot-z}  \\
		\dot{y} =-\varphi(y,t) \label{prot-y}
	\end{numcases}
\end{subequations}
with $y(0) = \nabla \mathcal{L}(z(0))$. Let $H = \nabla^2 F(x)$. The inversion of Hessian $\nabla\mathcal{L}(z)$ can be expanded as 
\begin{align}
	(\nabla^{2} \mathcal{L}(z))^{-1} &= \left[
	\begin{array}{cc}
		H & A^T\\
		A &  0 \\
	\end{array}
	\right]^{-1} \nonumber \\
	&= \left[
	\begin{array}{cc}
		H^{-1}-	H^{-1}A^TS^{-1}AH^{-1} &H^{-1}A^TS^{-1}\\
		S^{-1}AH^{-1} &  -S^{-1} \nonumber \\
	\end{array}
	\right]\\
	&= \left[
	\begin{array}{cc}
		H^{-1}-	H^{-1}A^TQ &Q^T\\
		Q &  -S^{-1} \label{Inversion-cen}\\
	\end{array}
	\right]
\end{align}
with $Q=S^{-1}AH^{-1} $ and $S= AH^{-1}A^T$ being the Schur complement of $H$ in $\nabla^{2} \mathcal{L}(z)$.
With properly chosen function $\varphi$, we can obtain the following convergence result.
\begin{ass}\label{ass-phi}The (Filippov) solution to subsystem \eqref{prot-y} (or \eqref{EZGS-y}) exists and \eqref{prot-y}  (or \eqref{EZGS-y}) is FT (resp. FxT/PT) stable at the origin for all solutions. 
\end{ass}
\begin{prop}\label{prop1}
Suppose that $F(x)$ is $\theta$-strongly convex and $A$ has full row rank. Under the dynamics \eqref{prot-cent}, for the solution $z(t)$ starting from $z_0$, it holds that  
\begin{enumerate}
	\item when the subsystem \eqref{prot-y} is exponentially stable with $\varphi(y)=r_0y, r_0>0$, there exists $\tau_0>0$ (dependent on $z_0$) such that $\|(\nabla^{2} \mathcal{L}(z(t)))^{-1}\| \leq \tau_0$ and $\|z(t)-z^* \|	\leq  \tau_0 y(0)e^{-r_0t}$; 
	\item when Assumption \ref{ass-phi} holds, $z(t)$ will converge to $z^*$ with the same FT (resp. FxT/PT), and $\|(\nabla^2\mathcal{L}(z(t)))^{-1}\|$ is bounded. 
\end{enumerate}
\end{prop}
\begin{proof}
	See Appendix\ref{App-CA}.
\end{proof}

To deal with the inequality constraints, we introduce the barrier functions to cast the inequality
constraints into the objective function. Motivated by \cite{fazlyab2017pred} or \cite[Sec. 11]{boyd2004convex}, we use the following approximate objective function:
\begin{align}\label{OB}
	\widetilde{F}(x,c(t),s(t)) \triangleq F(x)- \frac{1}{c(t)}\sum_{l=1}^{p}\log(s(t)-g^l(x)), 
\end{align} 
where $g^l(x)$ is the $l$-th component of $g(x)$, $c(t)>0$ a time-dependent positive barrier parameter, and $s(t) \geq 0$ is a time-dependent slack function used to relax the inequality constraints. Here, $c(t)$ and $s(t)$ are required to be differentiable over the time. For \eqref{OB}, the time-varying feasible set for the introduced barrier function is $\widetilde{D}(t) := \{x\in \mathbb{R}^n: g^l(x) < s(t), l \in \langle p \rangle \}$. Practically, one can choose $s(0)> \max_{l \in \langle p \rangle} g^{l}(x(0))$ such that $x(0) \in \widetilde{D}(0)$ for any initial state $x(0)$. Since $F(x)$ is strongly convex and the barrier functions are convex, then $\widetilde{F}(x,c(t),s(t))$ is strongly convex w.r.t. $x$. As the logarithmic barrier grows without bound if $s(t)-g^l(x) \rightarrow 0$, the optimal solution that minimizes \eqref{OB} subject to $Ax=b$ always exists and is unique, denoted by 
\begin{align}
	\tilde{x}^*(t) = {\arg\min}_{x \in \mathcal{S} \cap \widetilde{D}(t)}\widetilde{F}(x,c(t),s(t))
\end{align}
with $\mathcal{S} = \{x\in \mathbb{R}^n: Ax=b\}$. Introduce a modified Lagrangian function $\widetilde{\mathcal{L}}(z,c(t),s(t))=\widetilde{F}(x,c(t),s(t))+\lambda^T(Ax-b)$ and let $\tilde{z}^*(t)=(\tilde{x}^*(t), \tilde{\lambda}^*(t))$ be the corresponding saddle point. In order to track $\tilde{z}^*(t)$, the following dynamics with $z=(x,\lambda)$ is proposed 
\begin{subequations}\label{prot-cent2}
	%\small 
	\begin{numcases}{}
		\dot{z} =-\nabla_{zz}^{-1}\widetilde{\mathcal{L}} [\varphi(y,t)+\nabla_{zc} \widetilde{\mathcal{L}}\dot{c}+\nabla_{zs} \widetilde{\mathcal{L}}\dot{s} ], \ \ \quad \label{prot-z2}  \\
		\dot{y} =-\varphi(y,t) \label{prot-y2}
	\end{numcases}
\end{subequations}
with $y(0)=\nabla_z \widetilde{\mathcal{L}}(z(0),c(0),s(0))$. With properly chosen function $\varphi$, one can show that $z(t)$ tracks $\tilde{z}^*(t)$ with expected convergence performance, as indicated in the following result.
\begin{prop}\label{prop2}
	Suppose that $F(x)$ is $\theta$-strongly convex, $A$ has full row rank and Slater’s Condition holds. Consider the primal state $x(t)$ generated by \eqref{prot-cent2} with $x(0) \in \widetilde{D}(0)$ and $s(0)> \max_{l \in \langle p \rangle} g^{l}(x(0))$. When Assumption \ref{ass-phi} holds, it satisfies that $x(t) \in \widetilde{D}(t), \forall t\geq 0$ and $z(t)$ tracks $\tilde{z}^*(t)$ with the same FT (resp. FxT/PT) $T_0$. Moreover, when $s(t)$ converges to zero in FT (resp. FxT/PT) $T_1$ and $c(t) = c_0e^{r_ct}$, $x(t)$ will converge to $x^*$ exponentially with rate $r_c/2, \forall t\geq \max\{T_0, T_1\}$. 
\end{prop}
\begin{proof}
	See Appendix\ref{App-CB}.
\end{proof}

\section{Distributed EZGS Algorithms}\label{Sec-EZGS}
In this section, based on the previous Newton-based CTA \eqref{prot-cent}, we first provide an EZGS approach for designing distributed CTAs to solve \eqref{CO} with only equality constraints. Then, the convergence analysis is provided and the specific EXP/FT/FxT/PT convergence rates are shown for different EZGS algorithms, respectively. Finally, the algorithms are modified by incorporating local inequality constraints based on the log-barrier function.

\subsection{The EZGS Approach}
Consider the DO \eqref{CO} without inequality constraints, for which the Lagrangian function is given as 
\begin{align}
	\mathcal{L}(x,\lambda) = \sum_{i \in \mathcal{V} } f_i(x) + \sum_{i \in \mathcal{V} } \lambda_i^T(A_ix-b_i),
\end{align}
where $\lambda =[\lambda_i]_{i \in \mathcal{V}} \in \mathbb{R}^r$ and $\lambda_i$ is the Lagrangian multiplier corresponding to $A_ix=b_i$. By KKT conditions, $(x^*,\lambda^*)$ with $\lambda^* =[\lambda_i^*]_{i \in \mathcal{V}}$ is the optimal primal-dual solution iff the following conditions are satisfied 
\begin{align}\label{KKT-eq}
	\sum_{i \in \mathcal{V} } (\nabla f_i(x^*) + A_i^T \lambda_i^*) =0;
	A_ix^*-b_i = 0, \ \forall i \in \mathcal{V}. 
\end{align}
Let $x_i \in \mathbb{R}^n$ be the local copy of the global variable $x$ at agent $i$. Denote $\bm{z}=(\bm{x},\lambda)\in \mathbb{R}^{Nn+r}$ with $\bm{x}=[x_i]_{i\in \mathcal{V}} \in \mathbb{R}^{Nn}$ and $\overline{\mathcal{Z}}^*=\{\bm{1}_N \otimes x^*\}\times \mathcal{D}^*$. Moreover, when $A$ has full row rank, there exists a unique optimal primal-dual point represented by $z^*=(x^*,\lambda^*)$. Introduce the local Lagrangian function as $\mathcal{L}_i(x_i,\lambda_i) =  f_i(x_i) + \lambda_i^T(A_ix_i-b_i)$ and let $\widehat{\mathcal{L}}(\bm{x},\lambda)=\sum_{i=1}^{N} \mathcal{L}_i(x_i,\lambda_i)$. Partially inspired by the ZSG method \cite{Lu2012TAC}, we introduce three sets respectively defined as  
\begin{align}
	\mathcal{A} &=\{\bm{x} | x_1=x_2=\dots=x_N\},\\
	\mathcal{M} &=\{(\bm{x},\lambda) | \sum_{i \in \mathcal{V} } (\nabla f_i(x_i) + A_i^T \lambda_i) =0\}, \label{Set-M}\\
	\mathcal{C} &=\{\bm{x} | A_ix_i-b_i = 0, \forall i \in \mathcal{V} \} ,\label{Set-C}
\end{align}
where $\mathcal{A}$ represents the \textit{agreement set}, $\mathcal{M}$ denotes the \textit{ZGS manifold} for $\mathcal{L}_i$ w.r.t. $x_i$, i.e., $\sum_{i=1}^N \nabla_{x_i} \mathcal{L}_i(x_i,\lambda_i)=0, \forall (\bm{x},\lambda)\in \mathcal{M}$, and $\mathcal{C}$ represents the \textit{feasible constraint set}. According to \eqref{KKT-eq}, for the point $\bm{z}=(\bm{x},\lambda)$, if $\bm{z}\in \mathcal{M}$ and $\bm{x}\in \mathcal{A} \cap \mathcal{C}$, then $\bm{x}=\bm{1}_N \otimes x^*$ and $\lambda= \lambda^*$.
% where $(x^*,\lambda^*)$ is the optimal primal-dual pair to \eqref{CO}. 

Let $z_i=(x_i,\lambda_i)$ be the local state of agent $i$. To solve \eqref{CO} in a distributed way, the following dynamics is provided
\begin{subequations}\label{EZGS}\small
\begin{align}
	\dot{z}_i(t) &= -(\nabla^2 \mathcal{L}_i(z_i))^{-1} (g_i(y_i,t) + k(t)\sum_{j\in \mathcal{N}_i}\phi_{ij}(z_i,z_j)),\label{EZGS-z}\\
	\dot{y}_i(t) &= -g_i(y_i,t), \ \forall i\in \mathcal{V}  \label{EZGS-y}
\end{align}
\end{subequations}
with initialization $z_i(0)=(x_i(0),\lambda_i(0)) \in \mathbb{R}^n\times \mathbb{R}^{m_i}$, $y_i(0) = \nabla  \mathcal{L}_i(z_i(0)) $ and the time-varying gain function $k(t)$ to be determined.The right-hand side of \eqref{EZGS} is regarded as the control law $u_i(t)$ of the agent $i$. Denote $g_i(y_i,t)=[g_{x_i}(y_i,t);g_{\lambda_i}(y_i,t)]$ and $y_i = [y_{x_i};y_{\lambda_i}]$ to match the state $(x_i, \lambda_i)$. The local function $g_i$ can be chosen such that \eqref{EZGS-y} is EXP/FT/FxT/PT stable at the origin, e.g., $g_i(y) = a_i \text{sgn}^{\alpha_i}(y)+b_i\text{sgn}^{\beta_i}(y)$ with properly chosen coefficients $a_i,b_i, \alpha_i$ and $\beta_i$. Besides, $k(t)$ is a positive time-dependent function to be determined, and the local coupling function $\phi_{ij}(x_i,x_j)\triangleq [\chi_{ij}(x_i,x_j);\bm{0}_{r_i}]$ satisfying 
for any $(i,j) \in \mathcal{E}$
\begin{align}
	&\chi_{ij}(x_i,x_j)=-\chi_{ji}(x_j,x_i), \\
	&(x_i-x_j)^T\eta_{ij} >0, \forall \eta_{ij}\in \mathcal{F}(\chi_{ij})(x_i,x_j), \forall x_i\neq x_j. 
\end{align} 
Here, $\chi_{ij}$ is supposed to be measurable and locally essentially bounded to cover a wide range of functions, and $\mathcal{F}(\chi_{ij})$ denotes its Filippov set-valued map. Note that when $\chi_{ij}$ is continuous, we simply have $\mathcal{F}(\chi_{ij})(x_i,x_j)=\chi_{ij}(x_i,x_j)$. For example, one can choose $\chi_{ij}(x_i,x_j)=a_{ij}( \text{sgn}^{\alpha_{ij}}(e_{ij})+\eta\text{sgn}^{\beta_{ij}}(e_{ij}))$ with $e_{ij}=x_i-x_j$ and properly chosen coefficients $\eta, \alpha_{ij}$ and $\beta_{ij}$. With Assumptions \ref{ass-f} and \ref{ass-A}, $\nabla^2 \mathcal{L}_i(z_i)$ is invertible and \eqref{EZGS} is well defined for some properly chosen $k(t)$.  As only local and neighboring states are used in \eqref{EZGS}, it can be implemented in a distributed mode. 
%With the conditions on $g_i$ and $\chi_{ij}$, the Filippov solution always exists to \eqref{EZGS}.

In this work, we call the well-defined dynamics \eqref{EZGS} with $y_i(0) = \nabla  \mathcal{L}_i(z_i(0)) $ and aforementioned conditions on functions $g_i, \phi_{ij}, k(t)$ as an \textit{Extended Zero-Gradient-Sum} (EZGS) algorithm to be distinct from the ZGS method provided in \cite{Lu2012TAC} with $x_i(0)=x_i^*=\text{argmin}_x f_i(x)$, targeted for solving unconstrained DO. Differently from \cite{Lu2012TAC}, the proposed EZGS algorithm has free initial state $z_i(0)$ and can be used to solve a general class of equality constrained DO. Moreover, it can be extended to solve \eqref{CO} with general constraints by using barrier method, which will be discussed in Subsection \ref{Sec-DO-ineq}. Note that the protocol \eqref{EZGS} can be extended to the case where only a subset of $\mathcal{V}$ has equation constraints, and for those without equation constraints, the components corresponding to dual states are absent.   
Specially, when there are no equality constraints for all agents, with $x_i(0)=x_i^*$ indicating $y_i(0)=\nabla f_i(x_i(0))=0$, \eqref{EZGS} will reduce to ZGS methods considered in \cite{Lu2012TAC,song2016finite,Wu2021}. For the convenience of analysis, we suppose that all the agents have constraints.

\begin{ass}\label{ass-G}
	The undirected graph $\mathcal{G}(\mathcal{V}, \mathcal{E})$ is connected. 
\end{ass}
\subsection{Convergence Analysis}
In this subsection, we will show that the EZGS algorithm \eqref{EZGS} converges to $\mathcal{Z}^*$ over a graph $\mathcal{G}$ satisfying Assumption \ref{ass-G}. When the subsystem \eqref{EZGS-y} is FT stable, all $y_i(t)$ will converge to zero in a finite time $T_0<+\infty$. As a result, $g_i(y(t))=0$ for $t\geq T_0$. Then, \eqref{EZGS-z} will reduce to the following dynamics
\begin{align}\label{EZGS-sub}
	\dot{z}_i(t) &= -(\nabla^2 \mathcal{L}_i(z_i))^{-1} (\sum_{j\in \mathcal{N}_i}\phi_{ij}(z_i,z_j)). 
\end{align}
Denote $H_i=\nabla^2 f_i(x_i)$. Similar to \eqref{Inversion-cen}, the inversion of local Hessian $(\nabla^2 \mathcal{L}_i(z_i))^{-1}$ can be explicitly expressed as 
\begin{align}
	(\nabla^2 \mathcal{L}_i(z_i))^{-1}
%	&= \left[
%	\begin{array}{cc}
%		H_i & A_i^T\\
%		A_i &  0 \\
%	\end{array}
%	\right]^{-1} \nonumber \\
%	&= \left[
%	\begin{array}{cc}
%		H_i^{-1}-	H_i^{-1}A_i^TS_i^{-1}A_iH_i^{-1} &H_i^{-1}A_i^TS_i^{-1}\\
%		S_i^{-1}A_iH_i^{-1} &  -S_i^{-1} \nonumber \\
%	\end{array}
%	\right]\\
	= \left[
	\begin{array}{cc}
		H_i^{-1}-	H_i^{-1}A_i^TQ_i &Q_i^T\\
		Q_i &  -S_i^{-1} \label{Inversion}\\
	\end{array}
	\right]
\end{align}
with $Q_i=S_i^{-1}A_iH_i^{-1} $ and $S_i= A_iH_i^{-1}A_i^T$ being the Schur complement of $H_i$ in $\nabla^2 \mathcal{L}_i(z_i)$. Then, one can obtain the dynamics of $x_i$ over $[T,\infty)$ as follows
\begin{align}\label{ZGS-x}
	\dot{x}_i(t) &= -H_i^{-1}(I-A_i^TQ_i)\sum_{j\in \mathcal{N}_i}\chi_{ij}(x_i,x_j). 
\end{align}
For \eqref{ZGS-x}, in can be verified that $\mathcal{S}_i$ is invariant since $A_i\dot{x}_i(t)=0$. Moreover, $A_i^TQ_i$ is idempotent and we have $\|I-A_i^TQ_i\|\leq 1$.
Theorem \ref{thm-ZGS} gives the general convergence analysis of the proposed EZGS algorithm \eqref{EZGS}, and meanwhile indicates that \eqref{EZGS} enjoys the ZGS property. 

\begin{theorem}\label{thm-ZGS}
Suppose that Assumptions \ref{ass-f}-\ref{ass-G} hold and the subsystem \eqref{EZGS-y} is FT stable with settling time bounded by $T_0$. Then, for the EZGS algorithm \eqref{EZGS} with continuous function $k(t)>k_0$ over $[0 \ \infty)$ for some constant $k_0>0$, it holds: (i) $\bm{z}(t) \in \mathcal{S}_0=\mathcal{M} \cap (\mathcal{C}\times \mathbb{R}^r), \forall t\geq T_0$; (ii) $\lim_{t \rightarrow \infty } x_i = x^*, \forall i\in \mathcal{V}$.
\end{theorem}
\begin{proof}
See Appendix\ref{App-A}.
\end{proof}

%From Theorem \ref{thm-ZGS}, the following special case can be implied directly for the system \eqref{EZGS-sub}. 
%\begin{corollary}
%Let $z_i^*$ be the saddle point of $\mathcal{L}_i(x_i,\lambda_i), \forall i \in \mathcal{V}$. Then, with the initialization $z_i(0)=z_i^*$, $\mathcal{S}_0=\mathcal{M} \cap (\mathcal{C}\times \mathbb{R}^r)$ is invariant for \eqref{EZGS-sub} and $\lim_{t \rightarrow \infty } x_i(t) = x^*, \forall i\in \mathcal{V}$.
%\end{corollary} 

\subsection{Convergence Rate}
Before giving the main results, we present some representations and lemmas.
Let $\overline{A}=\text{blkdiag}(A_1,\cdots,A_N)$, $\overline{P}(\bm{x})=\text{blkdiag}(P_1(x_i),\cdots,P_N(x_N))$ with $P_i(x_i) =H_i^{-1}(x_i)-	H_i^{-1}(x_i)A_i^TQ_i(x_i)$ from \eqref{Inversion}, $\overline{B} = B\otimes I_n$ with $B$ being the incidence matrix of $\mathcal{G}$, and $M(\bm{x}) = \overline{B}^T \overline{P}(\bm{x}) \overline{B}$. It is obvious that $\overline{P}(\bm{x})$ is positive semi-definite. In the later analyses, we omit $\bm{x}$ for $\overline{P}(\bm{x})$ and $M(\bm{x})$ when no confusion arises. Then, Lemma \ref{lem-opt} can be obtained as an extension of \cite[Lem. 4]{Zhou2019} when considering a general projection matrix. Specially, when $H_i=I$, Lemma \ref{lem-opt} will reduce to \cite[Lem. 4]{Zhou2019}. Here, the proof of Lemma \ref{lem-opt} is omitted as it is analogous to that of \cite[Lem. 4]{Zhou2019}. Note that the full row rank assumption is imposed on $A$ rather than just $A_i$ as in Assumption \ref{ass-A}. The latter is not sufficient to guarantee Lemma \ref{lem-opt}. For example, if $A_i=I, \forall i\in \mathcal{V}$, then $\overline{P} =0 $, which gives $\mathcal{I}(\overline{B}) \cap \mathcal{N}(\overline{P})=\mathcal{I}(\overline{B}) $. However, Assumption \ref{ass-A} is sufficient for the asymptotic convergence of \eqref{EZGS} as shown in Theorem \ref{thm-ZGS}.
\begin{lem}\label{lem-opt}
	Suppose that Assumptions \ref{ass-f}, \ref{ass-G} hold and the matrix $A=[A_i]_{i\in \mathcal{V}}$ has full row rank. Then, $\mathcal{I}(\overline{B}) \cap \mathcal{N}(\overline{P})=\{0\}$. 
\end{lem}

\begin{lem}\cite[Lem. 3]{Shi2020auto}\label{lem-lam}
	Let $f: \mathbb{R}^n\mapsto \mathbb{R}^n$ be a continuous sign-preserving function, i.e., $\text{sign}(f(w))=\text{sign}(w), \forall w\in \mathbb{R}^n$. Given any nonzero positive semidefinite matrix $M \in \mathbb{R}^{n\times n}$, it holds that 
	\begin{align}\label{eq-lambdaf}
		\lambda_{f}(M) \triangleq \inf_{0\neq w \perp \mathcal{N}(M)} \frac{f^T(w)Mf(w)}{f^T(w)f(w)}>0.
	\end{align}
Particularly, when $f(w)=w$, $\lambda_{f}(M)=\lambda_{2}(M)$, which is the smallest positive eigenvalue of $M$.  
\end{lem}

\begin{lem}\label{lem-lam2}
Let $M(\bm{x}) = \overline{B}^T \overline{P}(\bm{x}) \overline{B}$ and $f: \mathbb{R}^{nm}\mapsto \mathbb{R}^{nm}$ be a continuous sign-preserving function. With the same assumptions in Lemmas \ref{lem-opt}, $\lambda_{f}(M(\bm{x}))$ is a positive continuous function of $\bm{x}$. 
\end{lem}
\begin{proof}
By Lemma \ref{lem-opt}, it is obvious that $M(\bm{x}) \neq 0$ and $\mathcal{N}(M(\bm{x}))=\mathcal{N}(\overline{B})$. Introduce the set $\mathcal{U} = \{u= \frac{f(w)}{\|f(w)\|}: w\in \mathcal{I}(\overline{B}^T), w\neq 0\}$. Then we can obtain the following parametric optimization from \eqref{eq-lambdaf} by considering $\bm{x}$ as a parameter
\begin{align}\label{PO}
	\lambda_{f}(M(\bm{x})) = \min_{u \in \overline{\mathcal{U}}} u^TM(\bm{x})u, 
\end{align}
in which $\overline{\mathcal{U}}$ is the closed closure of $\mathcal{U}$. As $\overline{\mathcal{U}}$ is compact and the objective function is continuous, then the minimizer of the parametric optimization can be attained. Moreover,  \eqref{eq-lambdaf} and \eqref{PO} share the same minimization value when $M=M(\bm{x})$. By Berge’s maximum theorem \cite[Sec. 17.5]{guide2006infinite}, since $\overline{\mathcal{U}}$ is a constant compact set, the minimization value $\lambda_{f}(M(\bm{x}))$ is continuous with respect to $\bm{x}$. 
\end{proof}

With Lemmas \ref{lem-opt}-\ref{lem-lam2}, the detailed convergence rate analyses are shown in the following result by choosing different types of local functions, i.e., $g_i,\chi_{ij}$. 
\begin{theorem}\label{thm2}
Suppose that Assumptions \ref{ass-f}, \ref{ass-G} hold and $A$ has full row rank. Consider the proposed EZGS \eqref{EZGS}, where $k(t)\equiv k_0>0$, $g_i$ and $\chi_{ij}$ respectively take the following form 
\begin{align*}
	g_i(y_i) &= a_i (\text{sgn}^{\alpha_i}(y_i)+\eta \text{sgn}^{\beta_i}(y_i)), \\
	\chi_{ij}(x_i,x_j) &= a_{ij}( \text{sgn}^{\alpha_{ij}}(x_i-x_j)+\eta \text{sgn}^{\beta_{ij}}(x_i-x_j))
\end{align*}
with $a_i>0, \eta\in\{0,1\}, \alpha_i, \alpha_{ij}=\alpha_{ji}\in [0 \ 1]$, $\beta_i, \beta_{ij}=\beta_{ji}>1$ and $a_{ij}$ is the weight of the edge $(i,j)$. Then,
\begin{enumerate}
	\item if $\eta=0$ and $\alpha_i=\alpha_{ij}=1$, $\bm{z}(t)$ will converge to $\overline{\mathcal{Z}}^*$ exponentially; 
	\item if $\eta=0$ and $\alpha_i\in [0 \ 1), \alpha_{ij}=0$, $\bm{z}(t)$ will converge to $\overline{\mathcal{Z}}^*$ in finite time; 
	\item if $\eta\in \{0,1\}$ and $\alpha_i\in [0 \ 1), \alpha_{ij}\in (0 \ 1)$, $\bm{z}(t)$ will converge to $\overline{\mathcal{Z}}^*$ in finite time;
	\item if $\eta=1$, $\alpha_i\in [0 \ 1), \alpha_{ij}\in (0 \ 1)$ and there exists a scalar $\underline{\lambda}_0$ such that $\lambda_{\chi}(M(\bm{x}))\geq \underline{\lambda}_0, \forall \bm{x}\in \mathbb{R}^{nN}$ with $\chi =[\chi_{ij}]_{(i,j)\in \mathcal{E}}$, $\bm{z}(t)$ will converge to $\overline{\mathcal{Z}}^*$ in fixed time. 
\end{enumerate}
%Moreover, if there exists a scalar $\underline{\lambda}_0$ such that $\lambda_{\chi}(M(\bm{x}))\geq \underline{\lambda}_0$ with $\chi =[\chi_{ij}]_{(i,j)\in \mathcal{E}}$ with $\eta=1$, $z(t)$ will converge to $z^*$ in fixed time. 
%Besides, if $g_i(y_i)=(c_i + \frac{ \dot{\mu}(t; t_0, T_0)}{ \mu(t; t_0, T_0)})y_i $ and $\chi_{ij}(x_i,x_j) = a_{ij}(d + \frac{ \dot{\mu}(t; t_0, T)}{ \mu(t; t_0, T)})(x_i-x_i)$ with $c_i, d>0$ and a prescribed time interval $T>T_0>t_0$, then $z(t)$ will converge to $\overline{\mathcal{Z}}^*$ in prescribed time $T$. 
\end{theorem}
\begin{proof}
	See Appendix\ref{App-B}.
\end{proof}

Theorem \ref{thm2} indicates that with different user-defined functions $g_i$ and $\chi_{ij}$, one can achieve EXP/FT/FxT/PT convergence. Moreover, the proposed FT/FxT EZGS algorithms are fully distributed with heterogeneous coefficients $\alpha_{ij}/\beta_{ij}$ to be distinct from the existing works \cite{song2016finite,Wu2021,Wang2020auto}. Note that the FxT convergence may not be ensured for the case 3) with $\eta=1$ as $\lambda_0$ (defined in the proof of Theorem \ref{thm2}) may depend on the initial states. For the case 4), as $\lambda_0$ is bounded by a positive constant, the FxT convergence can be achieved. For example, when $f_i$ is quadratic, $M$ is constant and the condition of case 4) holds directly. 

\begin{remark}
As equation constraints are considered, the Lyapunov candidate $V$ defined by \eqref{V0} is no longer used in the convergence rate analysis in Theorem \ref{thm2}, which is distinct from the existing works \cite{Lu2012TAC,song2016finite,Wu2021} for ZGS algorithms. The reason is that it is difficult to make $V(\bm{x})$ upper bounded by $\bm{x}^T\widetilde{L}\bm{x}$ where $\widetilde{L}$ denotes the Laplacian matrix for a weighted complete graph, as $x_i$ is constrained in $\mathcal{S}_i$. Moreover, we consider a more general case that the initial states are not required to be the local optima, which also covers the existing ZGS algorithms \cite{Lu2012TAC,song2016finite,Wu2021} when the constraints are absent and $x_i(0)=x_i^*$. 
Specifically, when there are no equations, the proposed FxT EZGS will reduce to \cite{Guo2023} with free initialization. 
\end{remark}

 Besides the FT/FxT convergence, the PT convergence rate is preferred in some real-world situations. For the completeness of the results, the PT convergent EZGS algorithm \eqref{EZGS} with typical control protocols will be presented. For this purpose, the following time-varying scaling function $\mu(t)$ is introduced \cite{wang2018prescribed}:  
\begin{equation}\label{mu}
	\mu(t;T)=\left\{\begin{aligned}
		\frac{T^h}{(T-t)^h}, \quad &t \in [0 \ T)\\
		1, \qquad \quad  & t \in [T \ \infty) 
	\end{aligned}\right.
\end{equation}
with $h>1$ and $T$ being the user-defined time. Besides, it holds that 
\begin{equation}\label{dmu}
	\dot{\mu}(t;T)=\left\{\begin{aligned}
		\frac{h}{T}\mu^{1+\frac{1}{h}}, \quad &t \in [0 \ T)\\
		0, \qquad \quad  & t \in [T \ \infty) 
	\end{aligned}\right.
\end{equation}
where the right-hand derivative of $\mu(t)$ is applied at $t=T$. For the analysis of PT ZGS algorithms, the following result is used, which is an extension of \cite[Lem. 1]{wang2018prescribed}, where only the undisturbed case $\bar{d}=0$ in \eqref{V-PT} is considered. 
\begin{lem}\label{lem-PT}
	Let $x(t; x_0)$ be any solution of the system (\ref{auto_dynamic}) starting from $x_0\in  \mathbb{R}^n$. Suppose that there exists a continuously differentiable function $V: \mathbb{R}^n\times \mathbb{R}_+ \rightarrow \mathbb{R}_+$ satisfying $V(0,t)=0$ and $V(x,t)>0, \forall x\neq 0$ such that 
	\begin{align}\label{V-PT}
		\dot{V} \leq -cV -2 \frac{\dot{\mu}(t;T)}{\mu(t;T)}V+\bar{d}\mu^l(t;T), 
	\end{align} 
	where $c\geq 0, \bar{d}\geq 0, l<\frac{1}{h}$, $\mu(t;T)$ and $\dot{\mu}(t;T)$ defined by \eqref{mu} and \eqref{dmu}, respectively. Then,  when $t\in [0, T)$, it holds that 
	\begin{enumerate}
	\item if $l = -2+\frac{1}{h}$
	\begin{align}\label{PT-Vderi1}
		V(t)  \leq \mu^{-2}e^{-ct}V(0)+\bar{d}T \mu^{-2} \text{ln}(\frac{T}{T-t}); 
	\end{align}
   \item if $l \neq -2+\frac{1}{h}$ 
   	\begin{align}\label{PT-Vderi11} 
   	V(t)  \leq \mu^{-2}e^{-ct}V(0)+\frac{\bar{d}T(\mu^{l-\frac{1}{h}}-\mu^{-2})}{(2+l)h-1}, 
   \end{align} 
    \end{enumerate}
where $\mu(t) =\mu(t;T)$. Moreover, $V(t)\equiv 0, \forall t\geq T$. Then, the origin of system \eqref{auto_dynamic} is globally PT stable with the prescribed time $T$. 
\end{lem}
\begin{proof}
For $t\in [0 \ T)$, multiplying $\mu^2$ on both sides of \eqref{V-PT} gives that 
\begin{align}
	\mu^2\dot{V} \leq -c\mu^2 V -2 \mu \dot{\mu} V+\bar{d}\mu^{2+l}, 
\end{align}
which implies that 
\begin{align}
\frac{d(\mu^2V)}{dt} \leq -c\mu^2 V +\bar{d}\mu^{2+l}. 
\end{align}
Solving the above differential inequality results in 
\begin{align*}
\mu^2(t)V(t) \leq e^{-ct} \mu^2(0)V(0) + \bar{d}\int_{0}^{t} e^{-c(t-s)} \mu^{2+l}(s)ds.
 %&\leq e^{-ct} V(0) + \bar{d}\int_{0}^{t}  \mu^2(s) ds \\
 %& =  e^{-ct} V(0) +  \frac{\bar{d}T}{(2+l)h-1}(\mu^{2+l-\frac{1}{h}}-1),
\end{align*}
Then, if $l = -2+\frac{1}{h}$, the above equation gives that 
\begin{align*}
	\mu^2(t)V(t) & \leq e^{-ct}V(0) + \bar{d}\int_{0}^{t} \frac{T}{T-t}(s)ds \\
	 & = e^{-ct}V(0) + \bar{d}T \text{ln}(\frac{T}{T-t}), 
	%&\leq e^{-ct} V(0) + \bar{d}\int_{0}^{t}  \mu^2(s) ds \\
	%& =  e^{-ct} V(0) +  \frac{\bar{d}T}{(2+l)h-1}(\mu^{2+l-\frac{1}{h}}-1),
\end{align*}
which yields \eqref{PT-Vderi1}. If $l \neq -2+\frac{1}{h}$, it holds that 
\begin{align*}
	\mu^2(t)V(t) \leq  e^{-ct} V(0) +  \frac{\bar{d}T}{(2+l)h-1}(\mu^{2+l-\frac{1}{h}}-1),
\end{align*}
which leads to \eqref{PT-Vderi11}. As $\lim_{t \rightarrow T^-} V(t) =0$ and $\dot{V}\leq 0$, we further have $V(t)=0, \forall t\geq T$. Consequently, the origin of system \eqref{auto_dynamic} is globally prescribed-time stable with the prescribed time $T$ as $l < \frac{1}{h}$. 
\end{proof}

Based on Lemma \ref{lem-PT}, the PT convergence of the EZGS algorithm with typical protocol is established in the following result. 
\begin{theorem}\label{thm3}
	Suppose that Assumptions \ref{ass-f}, \ref{ass-G} hold and $A$ has full row rank. Consider the proposed EZGS \eqref{EZGS} with initialization $\bm{x}_0$. Let $g_i(y_i)=(d + \frac{ \dot{\mu}(t; T_0)}{ \mu(t; T_0)})y_i $ and $\chi_{ij}(x_i,x_j) = a_{ij}(d + \kappa \frac{ \dot{\mu}(t; T)}{ \mu(t; T)})(x_i-x_i)$ with $c_i, d>0, \kappa >0 $ and a PT interval $T\geq T_0>0$. Then, there exists a positive scalar $\lambda_0$ (possibly dependent  on $\bm{x}_0$) such that when $\kappa \geq 1/\lambda_0$, $\bm{z}(t)$ will converge to $\mathcal{M}$ in PT $T_0$ and then converge to $\overline{\mathcal{Z}}^*$ in PT $T$. Besides, the control input of EZGS \eqref{EZGS} is continuous and uniformly bounded over $[0, \infty)$. 
\end{theorem}
\begin{proof}
	See Appendix\ref{App-C}.
\end{proof}

\begin{remark}
Compared with the EXP/FT/FxT protocols considered in Theorem \ref{thm2}, a time-varying gain $\kappa \frac{ \dot{\mu}(t; T)}{ \mu(t; T)}=\frac{\kappa h}{T-t}$ that tends to infinity as $t$ approaches $T$ is used to achieve PT convergence. Nevertheless, the control input of EZGS \eqref{EZGS} is also bounded over the time, making it feasible for practical implementation. Moreover, the PT convergence can be achieved in one stage, i.e., $T_0=T$. To guarantee that the inputs is bounded and continuous, it is necessary to choose the gain $\kappa$ larger than some positive scalar. For example, in \cite[Remark 3]{Zhou2020AUTO}, it has been shown that for the simplest system 
$\dot{x} = ax/(T-t)$, the input is bounded and continuous iff $a>1$. Moreover, according to the proof of Theorem \ref{thm3}, when there exists a scalar $\underline{\lambda}_0$ such that $\lambda_{2}(M(\bm{x}))\geq \underline{\lambda}_0, \forall \bm{x}\in \mathbb{R}^{nN}$, e.g., $f_i$ is quadratic, $\lambda_0$ in Theorem \ref{thm3} can be chosen as $\underline{\lambda}_0$ independent of the initial states. 

\end{remark}
\subsection{Extension to DO with Inequalities}\label{Sec-DO-ineq}
In this section, we will extend the EZGS algorithm \eqref{EZGS} for solving DO \eqref{CO} with general constraints. Referring to \eqref{OB}, based on log-barrier functions, we incorporate the inequality constraints into local objective functions, resulting in the following approximate version:
\begin{align}\label{OB-i}
	f_i^{c}(x,s_i) \triangleq f_i(x)- \frac{1}{c}\sum_{l=1}^{p_i}\log(s_i(t)-g_i^l(x)), \forall  i\in \mathcal{V}
\end{align} 
where $c>0$ and $g_i^l$ represents the $l$-th component of $g_i$. Note that $f_i^{c}$ is strongly convex and twice continuously differentiable w.r.t $x$ under the assumptions on $f_i$ and $g_i$ in Section \ref{sec2-A}. For \eqref{OB-i}, the time-varying feasible set for the introduced barrier function is $\widetilde{D}_i(t) := \{x\in \mathbb{R}^n: g_i^l(x) < s_i(t), l \in \langle p_i \rangle \}$. Practically, one can choose $s_i(0)> \max_{l \in \langle p_i \rangle} g_i^{l}(x(0))$ such that $x(0) \in \widetilde{D}_i(0)$ for any initial state $x_i(0)$.
Then, we get an approximate smooth optimization of \eqref{CO} with equality constraints
\begin{align}\label{CO-app}
	\min_{x} \sum_{i=1}^{N} f_i^{c}(x,s_i), \ \text{s.t.} \ x \in \cap_{i\in \mathcal{V}} \mathcal{S}_i.
\end{align}
Let $x_i$ be the local copy of $x$ and $\mathcal{L}_i^{c}(x_i,\lambda_i,s_i(t)) =  f_i^{c}(x,s_i) + \lambda_i^T(A_ix_i-b_i)$ be the local Lagrangian function to \eqref{CO-app}. By replacing $f_i$ in \eqref{Set-M} with $f_i^{c}$, one can defined a new ZGS set $\mathcal{M}^c$. Let $x_c^*$ be the optimal solution to \eqref{CO-app} subject to the parameter $c>0$ and $s_i=0$. Based on the analysis in \cite[Sec. 11.2]{boyd2004convex}, one can derive that $F(x_c^*) - F(x^*) \leq \frac{p}{c}$. 
%\begin{align}\label{dual-gap}
%	F(x_c^*) - F(x^*) \leq \frac{p}{c} .
%\end{align}
Since $F$ is $\underline{\theta}$-strongly convex with $\underline{\theta}=\min_{i\in \mathcal{V}}\{\theta_i\}$, it indicates that $\|x_c^*(t)-x^*\|\leq \sqrt{\frac{2p}{\underline{\theta}c}}$.
%\begin{align*}
%	\|\tilde{x}^*(t)-x^*\|\leq \sqrt{\frac{2p}{\underline{\theta}c}}. 
%\end{align*}
To obtain an approximate solution $x_c^*$ and facilitate the analysis, we consider a constant penalty scalar $c>0$. Then, inspired by \eqref{prot-cent2} and \eqref{EZGS}, we provide the following distributed EZGS dynamics 
\begin{subequations}\label{EZGS2}
	\begin{align}
		\dot{z}_i(t) &= -(\nabla^2 \mathcal{L}_i(z_i))^{-1} \mathcal{L}_i^c[g_i(y_i) + \sum_{j\in \mathcal{N}_i}\phi_{ij}(z_i,z_j) \nonumber \\ 
		& \qquad \qquad \qquad \qquad\qquad \qquad  +\nabla_{z_is_i} \mathcal{L}_i^c\dot{s}_i)],\label{EZGS2-z}\\
		\dot{y}_i(t) &= -g_i(y_i), \ \forall i\in \mathcal{V}  \label{EZGS2-y}
	\end{align}
\end{subequations}
with $y_i(0)=\nabla_{z_i} \mathcal{L}_i^c(z_i(0),s_i(0)), \forall i \in \mathcal{V}$.
\begin{theorem}
	Suppose that Assumptions \ref{ass-f}, \ref{ass-G} hold and $A$ has full row rank. Consider the modified EZGS \eqref{EZGS2} with $x_i(0) \in \widetilde{D}_i(0)$ and $s_i(0)> \max_{l \in \langle p_i \rangle} g_i^{l}(x(0))$. When the subsystem \eqref{EZGS2-y} is stable, $x_i(t) \in \widetilde{D}(t), \forall t\geq 0$. Moreover, when \eqref{EZGS2-y} is FT (resp. FxT/PT) stable with settling time bounded by $T_1$, $z(t)$ will converge to $\mathcal{S}_0^c=\mathcal{M}^c \cap (\mathcal{C}\times \mathbb{R}^m)$ in FT (resp. FxT/PT) $T_1$. Beside, if $s_i(t)$ converges to zero in FT (resp. FxT/PT) bounded by $T_2$,  with the FT (resp. FxT/PT) protocol provided in Theorem \ref{thm2}/\ref{thm3}, $x(t)$ will converge to $x_c^*$ in some FT (resp. FxT/PT) $T_3 > \max\{T_1,T_2\}$. 
\end{theorem}
\begin{proof}
Note that $\nabla_{x_i} \mathcal{L}_i^c(z_i(t))$ can be expressed as 
\begin{align}\label{grad-Li}
\nabla_{x_i} \mathcal{L}_i^c(z_i) = \nabla f_i(x_i)+\frac{1}{c} \sum_{l=1}^{p_i}\frac{\nabla g_i^l(x_i)}{s_i(t)-g_i^l(x_i) }+A_i^T\lambda_i.
\end{align}
Similar to the proof of Theorem \ref{thm-ZGS}, it can be shown that 
	\begin{align}
	\sum_{i \in \mathcal{V} }\nabla_{x_i} \mathcal{L}_i^c(z_i(t))= \sum_{i \in \mathcal{V} }y_{x_i}(t);
	\nabla_{\lambda_i} \mathcal{L}_i^c(z_i(t))=y_{\lambda_i}(t). \label{eqlam2}
\end{align}
Since the subsystem \eqref{EZGS2-y} is stable, one can conclude that each $\nabla_{x_i} \mathcal{L}_i^c(z_i(t))$ is bounded. It indicates that $x_i(t)$ never escapes $\widetilde{D}_i(t)$ by \eqref{grad-Li}, i.e., $x_i(t) \in \widetilde{D}_i(t), \forall t\geq 0$. Moreover, when \eqref{EZGS2-y} is FT (resp. FxT/PT) stable with settling time bounded by $T_1$, the right-hand sides of \eqref{eqlam2} converge to zero in the same FT (resp. FxT/PT), which means that $z(t)$ converges to $\mathcal{S}_0^c=\mathcal{M}^c \cap (\mathcal{C}\times \mathbb{R}^m)$ in FT (resp. FxT/PT) $T_1$. Besides, if $s_i(t)$ converges to zero in FT (resp. FxT/PT) bounded by $T_2$, i.e.,  $\dot{s}_i(t)=0, \forall t\geq T_2$, we have $x_i\in \mathcal{S}_i \cap \{x_i| g_i(x_i)<0\}, \forall t\geq \max\{T_1,T_2\}$, meaning that $x_i$ is locally feasible. Then, \eqref{EZGS2} will reduce to EZGS algorithm \eqref{EZGS}, and the remaining statement can be shown by applying Theorem \ref{thm2}/\ref{thm3}. 
\end{proof}

%\vspace{-1cm}
\begin{table}[b]
	\caption{Settings for four EZGS protocols}
	\label{tcl_prot}
	\centering
	%\normalsize
	\small
	\begin{tabular}{|l|l|l|}
		\hline
		Protocol     &  $\chi_{ij}(x_i,x_j), (i,j) \in \mathcal{E}$ &  	$g_i(y_i)$\\
		\hline
		LP &   $c_0(x_i-x_j)$  &  $c_0 y_i $   \\
		\hline
		FTP &   $5\text{sgn}^{\alpha_{ij}}(x_i-x_j)$  & $5\text{sgn}^{\alpha_i}(y_i)$  \\
		\hline
		FxTP &  $5(\text{sgn}^{\alpha_{ij}}+\text{sgn}^{\beta_{ij}})(x_i-x_j)$    &  $5(\text{sgn}^{\alpha_i}+\text{sgn}^{\beta_i})(y_i)$  \\
		\hline
		PTP &  $ (5 + \kappa \frac{ \dot{\mu}(t; 1)}{ \mu(t; 1)})(x_i-x_i)$    &  $(5 + \frac{ \dot{\mu}(t; 0.5)}{ \mu(t; 0.5)})y_i$  \\
		\hline
	\end{tabular}
\end{table}

%\begin{figure}[t]
%	\begin{center}
%		\subfigure[]{
%			\includegraphics[width=0.42\linewidth]{case1.pdf}
%			\label{Fig-case1}
%		}
%		\subfigure[]{
%			\includegraphics[width=0.44\linewidth]{case2.pdf}
%			\label{Fig-case2}
%		}
%		\caption{Values of $E_x$ and $E_{\lambda}$ for case 1 (a) and case 2 (b) by applying four protocols.}
%	\end{center}
%\end{figure}

\begin{figure}[t]
	\centering
	{\includegraphics[width=.35\textwidth]{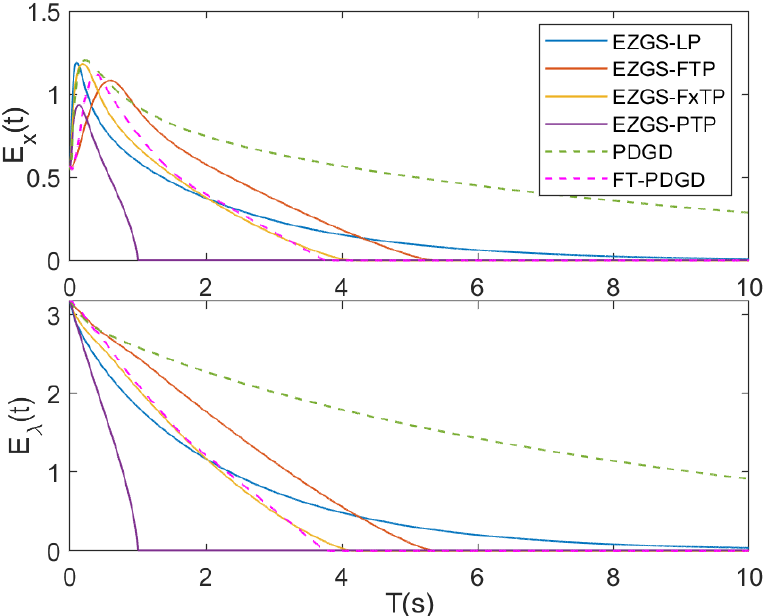}}
	\caption{Case 1: Comparison of $E_x$ and $E_{\lambda}$.}
	\label{Fig-case1}
\end{figure}
	
\section{Numerical examples}\label{Numberical}
Consider the equality constrained DO \eqref{CO} comprised of $N=6$ agents. For the $i$th agent, its local function takes the non-quadratic form $f_i(x) = \|x\|^2 - i*\bm{1}^Tx+\cos(w_i^T*x/2)$ with each entry of $w_i \in \mathbb{R}^7$ randomly chosen in $(0 \ 1)$, and equality constraint $A_ix=b_i$ corresponds to the partition given below
\begin{align*}
	A = \left[
	\begin{array}{ccccccc}
		0&1&2&3&3&-1&2\\ \hdashline[2pt/2pt] 
		1&0&2&-1&2&1&2\\  \hdashline[2pt/2pt] 
		0&1&2&0&-1&-1&0\\  \hdashline[2pt/2pt] 
		2&-1&2&1&-1&2&3\\  \hdashline[2pt/2pt] 
		1&1&3&0&2&3&0\\  \hdashline[2pt/2pt] 
		2&3&2&-1&0&-1&-1\\  
	\end{array}
	\right], b=\left[
	\begin{array}{c}
      -1\\ \hdashline[2pt/2pt] 
      2\\ \hdashline[2pt/2pt] 
      2\\ \hdashline[2pt/2pt] 
      2\\ \hdashline[2pt/2pt] 
      2\\ \hdashline[2pt/2pt] 
      3\\ 
	\end{array}
	\right].
\end{align*}
The optimal primal-dual solution is $z^*=(x^*,\lambda^*)$ with $x^* \approx (-0.100, 0.763, 0.506, -0.710, -0.490, 0.266, 0.546)$ and $\lambda^*\approx (7.838,	-6.301,	-12.907, 5.947,	4.553, 6.165)$. Let $\mathcal{G}$ be a connected circle graph with edge weights being all ones, i.e., $a_{ij}=1, \forall (i,j)\in \mathcal{E}$; otherwise, $a_{ij}=0$. For the first study, four types of EZGS protocols \eqref{EZGS} are designed for solving the above DO, i.e., linear protocol, FT protocol, FxT protocol and PT protocol, abbreviated as LP, FTP, FxTP and PTP, respectively. The detailed settings for all protocols are given in Table \ref{tcl_prot}, where $c_0=20$, $\kappa=10$ and $h=3$ for the function $\mu(\cdot)$. The heterogeneous coefficients of FTP/FxTP are chosen as $\alpha_i=0.1*i$, $\beta_i=1+0.1*i$, $\alpha_{ij}=\alpha_{ji} = 0.1*\min\{i,j\}$ and $\beta_{ij}=\beta_{ji} =1+0.1*\min\{i,j\}$. The initial primal-dual states $z_i(0)=(x_i(0), \lambda_i(0))$ are all zeros and $y_i(0) = \nabla_{z_i}  \mathcal{L}(z_i(0)) $. The proposed EZGS protocols are further compared to the traditional primal-dual gradient descent dynamics (PDGD) \cite{Cherukuri2018} and FT-PDGD \cite{Shi2023Tcyber} adapted to equation-constrained DO, which respectively take the form as \cite[Alg. (41) with $\rho=1$]{Shi2023Tcyber} with $g_{x/\lambda/\mu}(w) =5*w $ and $g_{x/\lambda/\mu}(w)=5*\text{sign}(w)$ in the simulations. Compared with EZGS methods, additional dual variables corresponding to the consensus constraints are updated and communicated in PDGD and FT-PDGD, resulting in a higher communication burden. However, as the inversion of Hession matrix is used in the proposed EZGS methods, the computation costs of which are higher than those of PDGD and FT-PDGD.

For the sake of performance comparison, we define $E_x(t)= \frac{1}{N}\sum_{i=1}^{N}\|x_i(t)-x^*\|$ and $E_{\lambda}(t)= \frac{1}{N}\sum_{i=1}^{N}\|\lambda_i(t)-\lambda_i^*\|$ as the system error function. By the simulation, the values of $E_x$ and $E_{\lambda}$ by applying six protocols are shown in Fig. \ref{Fig-case1}. It can be seen that all the local primal and dual states with EZGS-FTP/FxTP and FT-PDGD converge to $x^*$ and $\lambda_i^*$ in FT/FxT, while only exponential convergence can be obtained with EZGS-LP and PDGD. Notably, with EZGS-PTP, the convergence to $x^*$ and $\lambda_i^*$ is achieved in PT $T=1s$, as shown in Fig. \ref{Fig-case1}. These results verify Theorems \ref{thm2} and \ref{thm3}.

\begin{figure}[t]
	\centering
	{\includegraphics[width=.35\textwidth]{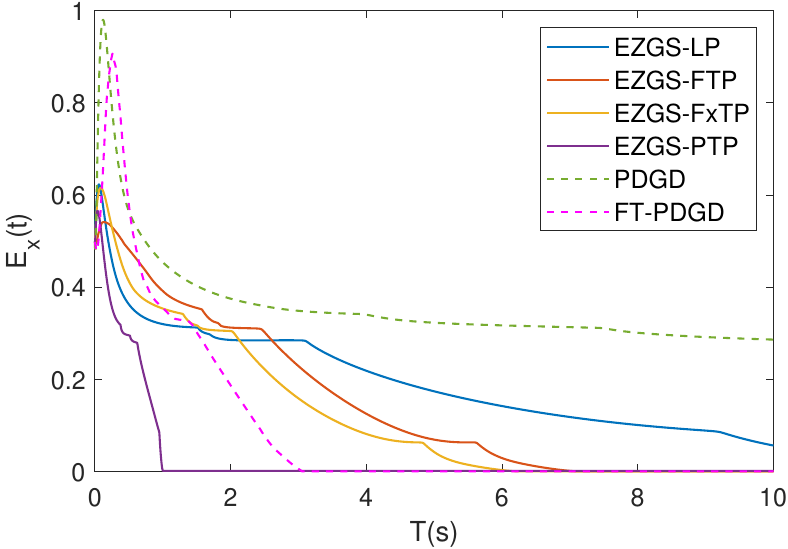}}
	\caption{Case 2: Comparison of $E_x$.}
	\label{Fig-case2}
\end{figure}

For the second case study, the local inequality constraint for the agent $i$ is added to the previous DO as $\bm{1}_N^Tx_i -x_{ii}\leq 1+(i-1)/10$, where $x_{ii}$ is the $i$-th entry of the local copy $x_i$. In this case, the optimal primal solution is $x^* \approx (-0.034, 0.540, 0.597, -0.683, -0.436, 0.169, 0.396)$. The modified EZGS \eqref{EZGS2} is applied with $c=1000, z_i(0)=\bm{0}$ and $s_i(t) \equiv 0$ as the initial states are feasible. Similar as the previous case, four types of EZGS protocols are used as shown in Table \ref{tcl_prot}, where $c_0=\kappa=20$ to enhance the convergence. For the sake of comparison, the PDGD and FT-PDGD algorithms with the same settings as the previous case are implemented. The simulation results are shown in Fig. \ref{Fig-case2}. It can be seen that states with FTP/FxTP/PTP converge to the approximate solution $x_c^*$ (near $x^*$) in FT/FxT/PT while only exponential convergence can be achieved with LP. Meanwhile, all the local primal states are feasible during the simulation. Note that PDGD and FT-PDGD algorithms can achieve asymptotic and FT convergence to $x^*$, respectively, but with more dual variables and communication cost.

\section{Conclusion}\label{conclusion}
In this paper, we have provided an EZGS approach for solving constrained DO with general constraints. We designed a centralized Newton-based CTA and several EZGS CTAs for solving centralized and distributed constrained optimization with free initialization, respectively, which can achieve the EXP/FT/FxT/PT convergence with typical protocols. By incorporating with the barrier method, the proposed CTAs have been modified to tackle general inequality constraints. Further work will focus on EZGS CTA over time-varying graphs and its discretization.

%In this paper, three FT/FxT distributed algorithms are developed over a multi-agent network for linear equations $Ax=b$ when $A$ has full row rank. With specific initializations, a nonlinear distributed protocol is provided by combining the projection-consensus flow \cite{Shi2017TAC} and a general version of the existing FT/FxT consensus algorithms \cite{Wang2010,zuo2018TII}, which can obtain a solution in FT/FxT time over an undirected and fixed network. To eliminate the initial conditions, a nonlinear term is added to drive the local states to satisfy the specific initialization. When multiple solutions emerge, a modified initialization is given for the predesigned distributed algorithms to obtain one solution closest to specific points in FT/FxT. For the purpose of practical applications, the command signals deigned from the continuous-time protocols can be converted to digital forms by using sampling technique. Moreover, the delay and asynchronous communication may occur in the network control systems. Future work will concern the design and stability analysis of robust protocols under the uncertain environment in presence of signal delay and noises as well as the asynchronous implementation of its discrete version for solving equation problems.

\appendices 

\section*{Appendix}
\subsection{Proof of Proposition \ref{prop1}}\label{App-CA}
	With \eqref{prot-cent}, we have $\dot{\nabla}\mathcal{L}(z(t)) ={\nabla^2\mathcal{L}}(z)\dot{z}(t)=-\varphi(y,t)= \dot{y}$. Since $y(0) = \nabla \mathcal{L}(z(0))$, there exists a solution $y(t)$ to \eqref{prot-y} such that $\nabla\mathcal{L}(z(t)) \equiv y(t) $. 
	
	(1) When the subsystem \eqref{prot-y} is exponentially stable with $\varphi(y)=r_0y, r_0>0$, we get $y(t) = y(0)e^{-r_0t}$. Applying the mean-value theorem to expand $\nabla \mathcal{L}(z)$ around $z^*=(x^*,\lambda^*)$, we have
	\begin{align}\label{mean0}
		\nabla \mathcal{L}(z) = \nabla^2\mathcal{L}(\eta(z))(z-z^*),
	\end{align}
	where $\eta(z)$ is a convex combination of $z$ and $z^*$. Then, it gives that
	\begin{align}\label{mean0-inv}
		\|z-z^* \|	\leq \|(\nabla^2\mathcal{L}(\eta(z)))^{-1} \|\|\nabla \mathcal{L}(z)\|.
	\end{align}
	Consider the following function 
	\begin{align}\label{V0-cen}
		V_0(x) =  F(x^*)-F(x)-\nabla^T F(x) (x^*-x),
	\end{align} 
	which satisfies $	V_0(x) \geq \frac{\theta}{2} \|x-x^*\|^2$ since $F$ is $\theta$-strongly convex. Let $e_x = x-x^*$ and $y=[y_x; y_{\lambda}]$. With \eqref{Inversion-cen}, taking the derivative of $V_0(x)$ along \eqref{prot-z} gives that 
	\begin{align*}
		\dot{V}_0 &=  -  (x^*-x)^T \nabla^2 F(x)\dot{x}\\
		&=- e_x^T(I-A^TQ)(r_0y_{x}(t)) - r_0(Ax-Ax^*)^T(S^{-1}y_{\lambda}(t))\\
		&=- e_x^T(I-A^TQ)(r_0y_x(t)) - r_0(Ax-b)^T(S^{-1}y_{\lambda}(t))\\
		&=- e_x^T(I-A^TQ)(r_0y_{x}(t)) - r_0(y_{\lambda})^T(S^{-1}y_{\lambda}(t))\\
		&\leq r_0e^{-r_0 t}\|y_{x}(0)\|\|e_x\| \leq r_0 \|y_{x}(0)\| e^{-r_0 t} \sqrt{2/{\theta}}\sqrt{V_0}.
	\end{align*} 
Then, one can further calculate that for all $t\geq 0$, $	V_0(t) \leq (\sqrt{V_0(x_0)}+\|y_x(0)\|/\sqrt{2\theta})^2$, 
%	\begin{align*}
%		V_0(t) \leq (\sqrt{V_0}+|y_x(0)|\sqrt{\frac{1}{2\theta}})^2,
%	\end{align*}
	which indicates that $x(t)$ is bounded, i.e., $x(t) \in U_0 =\{x: \|x-x^*\|\leq \sqrt{2V_0(x_0)/\theta}+y_x(0)/\theta\}$. As $(\nabla^{2}\mathcal{L}(z(t)))^{-1}$ is well-defined and continuous, there must exist a constant $\tau_0>0$ (dependent on $z_0$) such that $\|(\nabla^2\mathcal{L}(z(t)))^{-1}\| \leq \tau_0$ for any $x \in U_0$. As $x^*\in U_0$, with \eqref{mean0-inv}, we can obtain that 
	\begin{align*}
		\|z(t)-z^* \|	\leq  \tau_0 \|\nabla \mathcal{L}(z(t))\| =\tau_0 \|y(0)\|e^{-r_0t}.
	\end{align*}
	
	(2) When the subsystem \eqref{prot-y} is FT (resp. FxT/PT) stable at the origin, $\nabla\mathcal{L}(z)$ will converge to 0 with the same FT (resp. FxT/PT) $T_0$. It means that $z(t)$ will converge to $z^*$ at $T_0$ and stay thereby according to the optimal condition. As $z(t)$ is continuous over $[0, T_0]$, $\|(\nabla^2\mathcal{L}(z(t)))^{-1}\|$ is bounded over $[0, T_0]$. The proof is thus completed. 
	%Then, the boundedness of the right-hand side of \eqref{prot-cent} can be implied by that of the subsystem \eqref{prot-y} and $\varphi(y(t),t)=0, \forall t\geq T_0$.

\subsection{Proof of Proposition \ref{prop2}}\label{App-CB}
	By taking the time derivative of $\nabla_z \widetilde{\mathcal{L}}$ and referring to \eqref{prot-cent2}, we have 
	\begin{align*}
		\dot{\nabla}_z \widetilde{\mathcal{L}}=\dot{\nabla}_{zz} \widetilde{\mathcal{L}}\dot{z} + \nabla_{zc} \widetilde{\mathcal{L}}\dot{c}+\nabla_{zs} \widetilde{\mathcal{L}}\dot{s}
		= -\varphi(y,t) = \dot{y}. 
	\end{align*}
	Combining with $y(0) = \nabla \mathcal{L}(z(0))$, it holds that $\nabla_z\widetilde{\mathcal{L}} \equiv y(t) , \forall t\geq 0$ for some solution $y(t)$ to \eqref{prot-y2}. As the subsystem \eqref{prot-y2} is FT (resp. FxT/PT) stable at the origin with continuous term $\varphi(y(t),t)$, then $\nabla_z\widetilde{\mathcal{L}}$ will converge to 0 with the same FT (resp. FxT/PT) $T_0$, i.e., $\nabla_z\widetilde{\mathcal{L}}(z(t),c(t),s(t)) =0, \forall t\geq T_0$. It indicates that $z(t)=\tilde{z}^*(t)$ by the KKT condition when $t\geq T_0$ and $\|\nabla_z\widetilde{\mathcal{L}}(z(t),c(t),s(t))\|$ is bounded for all $t\geq 0$ as $\nabla_z\widetilde{\mathcal{L}}$ is continuous over $[0, T_0]$. As a result, $\|\nabla_x \widetilde{\mathcal{L}}(z(t),c(t),s(t))\|$ is bounded for $t\geq 0$, which implies that $x \in \widetilde{D}(t)$ as
	\begin{align}
		\nabla_x \widetilde{\mathcal{L}}(z(t),c(t),s(t)) = \nabla F(x) + \frac{1}{c}\sum_{l=1}^p \frac{\nabla g^l(x)}{s-g^l(x)}.
	\end{align}
	 Moreover, when $s(t)$ converges to zero in FT (resp. FxT/PT) $T_1$, we have $\widetilde{D}(t) := \{x\in \mathbb{R}^n: g^l(x) <0, l \in \langle p \rangle \}, \forall t \geq T_1$. Then, one can conclude that  $x(t)$ is feasible for any $t\geq T_2 \triangleq \max\{T_0,T_1\}$, i.e., $x(t)\in \mathcal{X}, \forall t \geq T_2$. According to \cite[Lemma 1]{fazlyab2017pred} with $s(t)=0, \forall t\geq T_2$, one can obtain that 
	\begin{align}
		F(\tilde{x}^*(t)) - F(x^*) \leq \frac{p}{c(t)}, \forall t \geq T_2.
	\end{align}
	Since $F$ is $\theta$-strongly convex, so is $F+\delta_{\mathcal{X}}$ with $\delta_{\mathcal{X}}$ being the indicator function for $\mathcal{X}$. Then, as $x^*$ minimizes $F+\delta_{\mathcal{X}}$, it holds that 
	\begin{align*}
		\frac{\theta}{2}\|\tilde{x}^*(t) -x^* \|^2&\leq F(\tilde{x}^*(t))+\delta_{\mathcal{X}}(\tilde{x}^*(t)) - F(x^*) -\delta_{\mathcal{X}}(x^*(t))\\
		&\leq \frac{p}{c(t)}= pc_0^{-1} e^{-r_ct}, \forall t \geq T_2.
	\end{align*}
	Combining with $x(t) = \tilde{x}^*(t), \forall t \geq T_2$, one can conclude that $x(t)$ will converge to $x^*$ exponentially with rate $r_c/2$.

\subsection{Proof of Theorem \ref{thm-ZGS}}\label{App-A}
	%Let $\widehat{\mathcal{L}}(\bm{x},\lambda)=\sum_{i=1}^{N} \mathcal{L}_i(x_i,\lambda_i)$. Then, 
	Since $\chi_{ij}$ is measurable, locally essentially bounded and $k(t)$ is continuous over $[0 \ \infty)$, with Assumption \ref{ass-phi}, the solution to the proposed EZGS \eqref{EZGS} always exists. Moreover, the derivative of $\nabla \mathcal{L}_i(z_i)$ is given by
	\begin{align*}
		\frac{d}{dt} \nabla \mathcal{L}_i(z_i)& = \nabla^2 \mathcal{L}_i(z_i)\dot{z}_i =- (g_i(y_i,t) + k(t)\sum_{j\in \mathcal{N}_i}\phi_{ij}(z_i,z_j)).
	\end{align*}
By the definition of Filippov solution and the fact $k(t)$ is continuous, there exist $\omega_{x_i}(y_i,t) \in \mathcal{F}[g_{x_i}](y_i,t)$ and $\eta_{ij}(x_i,x_j)(=-\eta_{ji}(x_j,x_i)) \in  \mathcal{F}[\chi_{ij}](x_i,x_j)$ such that 
	\begin{align*}
		\sum_{i \in \mathcal{V} }\frac{d}{dt} \nabla_{x_i} \mathcal{L}_i(z_i)& \overset{a.e.}= - \sum_{i \in \mathcal{V} }( \omega_{x_i}(y_i,t)  + k(t)\sum_{j\in \mathcal{N}_i}\eta_{ij}(x_i,x_j)), \\
		& \overset{a.e.}= - \sum_{i \in \mathcal{V} }\omega_{x_i}(y_i) = \sum_{i \in \mathcal{V} }\dot{y}_{x_i},
	\end{align*}
	and meanwhile $\frac{d}{dt} \nabla_{\lambda_i} \mathcal{L}_i(z_i) = \dot{y}_{\lambda_i}$ for some solution $y_i(t)$ to the subsystem \eqref{EZGS-y}. With the initialization $y_i(0)=\nabla_{z_i}  \mathcal{L}(z_i(0))$, we have $\sum_{i=1}^{N} \nabla_{x_i} \mathcal{L}_i(z_i(0))= \sum_{i=1}^{N}y_{x_i}(0)$ and $\nabla_{\lambda_i} \mathcal{L}_i(z_i(0))=y_{\lambda_i}(0)$. Then, it can be obtained that for any $t\geq 0$
	\begin{align}
		\sum_{i \in \mathcal{V} }\nabla_{x_i} \mathcal{L}_i(z_i(t))&= \sum_{i \in \mathcal{V} }y_{x_i}(t), \label{eqx}\\
		\nabla_{\lambda_i} \mathcal{L}_i(z_i(t))&=y_{\lambda_i}(t). \label{eqlam}
	\end{align}
	As the subsystem \eqref{EZGS-y} is FT stable within time $T_0$, $y_i(t)=0$ and $g_i(y_i(t))=0, \forall t\geq T_0$, which indicates that $\sum_{i=1}^{N}\nabla_{x_i} \mathcal{L}_i(z_i(t))=0$ and $\nabla_{\lambda_i} \mathcal{L}_i(z_i(t)), \forall t\geq T_0$. Statement (i) holds directly.

	To show statements (ii), we suppose that $t\geq T_0$ in the following analysis. Consider the following Lyapunov candidate 
	\begin{align}\label{V0}
		V(\bm{x}) = \sum_{i \in \mathcal{V} } f_i(x^*)-f_i(x_i)-\nabla^T f_i(x_i) (x^*-x_i).
	\end{align} 
	Since $f_i$ is $\theta_i$-strongly convex, it holds that 
	\begin{align}\label{V_lb0}
		V(\bm{x}) \geq \sum_{i \in \mathcal{V} } \frac{\theta_i}{2}\|x_i-x^*\|^2\geq \frac{\underline{\theta}}{2}\sum_{i \in \mathcal{V} } \|x_i-x^*\|^2
	\end{align}
	with $\underline{\theta} = \min_{i\in \mathcal{V}}\{\theta_i\}$. With $y_i(0)=0$, \eqref{EZGS-z}  becomes \eqref{EZGS-sub} for $t\geq T_0$. Taking the time derivative of $V$ over $[T_0 \ \infty)$ along the dynamics \eqref{ZGS-x} gives that
	\begin{align*}
		\dot{V}(\bm{x}) &= -\sum_{i \in \mathcal{V} }  (x^*-x_i)^T \nabla^2 f_i(x_i)\dot{x}_i\\
		& \overset{a.e.}= -k(t)\sum_{i \in \mathcal{V} }  (x^*-x_i)^T(I-A_i^TQ_i)\sum_{j\in \mathcal{N}_i}\eta_{ij}(x_i,x_j) \\
	%	& = -\sum_{i=1}^N (x^*-x_i)^T\sum_{j\in \mathcal{N}_i}\eta_{ij}(x_i,x_j) \\
		& \leq -\frac{k_0}{2}\sum_{i \in \mathcal{V} }\sum_{j\in \mathcal{N}_i}(x_j-x_i)^T\eta_{ji}(x_j,x_i)\leq 0, 
	\end{align*}	
	where $\eta_{ji}(x_j,x_i) \in \mathcal{F}(\chi_{ji})(x_j,x_i)$ and the third equation holds since $A_i(x^*-x_i)=0, \forall t\geq T_0$. Then, $V(\bm{x}(t)) \leq V(\bm{x}(T_0))$ and hence
	\begin{align}\label{Phi}
		\bm{x}(t) \in \Phi\triangleq 
		\{\bm{x}:\frac{\underline{\theta}}{2}\sum_{i \in \mathcal{V} } \|x_i-x^*\|^2 \leq  V(\bm{x}(T_0))\}, \forall t\geq T_0.
	\end{align} 
	Since $f_i$ is twice continuously differentiable, then there exists $\Theta_i >0$ such that $H_i(x_i)\leq \Theta_i I, \forall i\in \mathcal{V}$ for $\bm{x} \in \Phi$. Then, $V(\bm{x})$ satisfies 
	\begin{align}\label{V_bb}
		\sum_{i \in \mathcal{V} } \frac{\theta_i}{2}\|x_i-x^*\|^2  \leq 	V(\bm{x}) \leq \sum_{i \in \mathcal{V} } \frac{\Theta_i}{2}\|x_i-x^*\|^2.
	\end{align}
	Note that $\dot{V}(\bm{x}) = 0$ iff $\bm{x} \in  \bm{1}_N \otimes \mathbb{R}^n$. Recall that $\bm{x}(t) \in \mathcal{X}_0 \triangleq \mathcal{M}|_{\bm{x}} \cap \mathcal{C}$ with $\mathcal{M}|_{\bm{x}}$ be the projection on the space of $\bm{x}$. 
	By the optimal condition, $(\bm{1}_N \otimes \mathbb{R}^n) \cap \mathcal{X}_0=\{\bm{1}_N \otimes x^*\} $, which indicates that $\dot{V}(\bm{x}) <0$ if $\bm{x}\neq \bm{1}_N \otimes x^*$. Then, according to \cite[Proposition 3]{Bacciotti2006}, 
	$\lim_{t \rightarrow \infty } \bm{x} = \bm{1}_N \otimes x^*$, which concludes the statement (ii). 
	%The proof is complete.

\subsection{Proof of Theorem \ref{thm2}}\label{App-B}
	Denote $\bar{\alpha} =[\alpha_{ij}]_{(i,j)\in \mathcal{E}}\otimes \bm{1}_n, \bar{\beta} =[\beta_{ij}]_{(i,j)\in \mathcal{E}}\otimes \bm{1}_n$, and $\bm{\xi} =[\xi_{ij}]_{(i,j)\in \mathcal{E}}= \overline{B}^T\bm{x} $ with $\xi_{ij} = x_i-x_j$ when $(i,j)\in \mathcal{E}$. Let $\chi(\bm{\xi}) =[\chi_{ij}]_{(i,j)\in \mathcal{E}} = \text{sgn}^{[\bar{\alpha}]}(\bm{\xi}) + \eta \text{sgn}^{[\bar{\beta}]}(\bm{\xi})$. Consider the following function 
	\begin{align}\label{V-1}
		V_1 =\sum_{(i,j)\in \mathcal{E}}\big(\frac{\mathbf{1}_n^T|\xi_{ij}|^{\alpha_{ij}+1} }{\alpha_{ij}+1}+ \frac{\eta\mathbf{1}_n^T|\xi_{ij}|^{\beta_{ij}+1}}{\beta_{ij}+1}\big), 
		%\\
		%&=\sum_{l=1}^{\bar{l}}\big(\frac{\hat{a}_l}{\alpha_{l}+1}|z_l|^{\alpha_{l}+1} + \frac{\hat{b}_l}{\beta_{l}+1}|z_l|^{\beta_{l}+1}\big).
	\end{align}	
	which will be used in the later convergence analysis. 	
	
	(1) If $\eta=0$ and $\alpha_i = \alpha_{ij}=1$, $g_i(y_i)=a_iy_i$ and $\chi_{ij}(x_i,x_j)=a_{ij}(x_i-x_j)$.
	Then, one has $y_i(t) = y_i(0)e^(-a_it)$. With \eqref{eqlam}, it means that $A_ix_i-b_i=A_i(x_i-x^*)=y_i(0)|_{\lambda_i}e^{-a_it}$. Using \eqref{Inversion}, the dynamics of $x_i$ can be expressed by 
	\begin{align*}\small
		\dot{x}_i(t) &= -H_i^{-1}(I-A_i^TQ_i)(a_iy_{x_i}+ \sum_{j\in \mathcal{N}_i}\chi_{ij}(x_i,x_j)) -a_iQ_i^T y_{\lambda_i}. 
	\end{align*}
	 Let $e_i=x_i-x^*$ and $\bm{e}=[e_i]_{i\in \mathcal{V}}$. Consider the function $V(\bm{x})$ expressed by \eqref{V0}, which satisfies 
	\begin{align}\label{V_relax}
		\dot{V}(\bm{x})=& -\frac{1}{2}\sum_{i \in \mathcal{V}}\sum_{j\in \mathcal{N}_i}(x_j-x_i)^T\chi_{ji}(x_j,x_i)- \sum_{i \in \mathcal{V} } (x_i-x^*)^T \nonumber \\
		& [(I-A_i^TQ_i)(a_iy_{x_i})+a_iA_i^TS_i^{-1}y_{\lambda_i}] \nonumber\\
		\leq & - \sum_{i \in \mathcal{V} } e_i^T(I-A_i^TQ_i)(a_iy_{x_i}) - a_iy_{\lambda_i}^TS_i^{-1}y_{\lambda_i} \\
	%	\leq &  \sum_{i \in \mathcal{V} } a_ie^{-a_it}\|y_{x_i}(0)\|\|e_i\|- a_iy_{\lambda_i}^TS_i^{-1}y_{\lambda_i}\\	
		\leq &\sum_{i \in \mathcal{V} } a_ie^{-a_it}\|y_{x_i}(0)\|\|e_i\|\leq w e^{-\underline{a} t} \sqrt{2N/\underline{\theta}}\sqrt{V}, \nonumber
	\end{align}	
	where $w = \max_{i \in \mathcal{V}}\{a_i \|y_i(0)|_{x_i}\|\}$ and $\underline{a} =\min_{i \in \mathcal{V}}\{a_i \} $. Then, it can be calculated that $V$ satisfies for any $t\geq 0$
	\begin{align}
		V(\bm{x}(t)) \leq (\sqrt{V(0)}+\underline{a}w\sqrt{N/(2\underline{\theta}}))^2 \triangleq \varepsilon_1 .
	\end{align}  
	It can be concluded that $U= \{\bm{x}\in \mathbb{R}^{Nn}: 	V(\bm{x}) \leq \varepsilon_1\}$ and $U_i= \{x_i\in \mathbb{R}^n: f_i(x^*)-f_i(x_i)-\nabla^T f_i(x_i) (x^*-x_i) \leq \varepsilon_1\}$ are compact. Obviously, $x^* \in U_i, \forall i \in \mathcal{V}$. Since $f_i$ is twice continuously differentiable, then there exists $\tilde{\Theta}_i >0$ such that $H_i(x_i)\leq \tilde{\Theta}_i I$ for any $i\in \mathcal{V}$, which indicates that $\theta_i \leq \|H_i\|\leq \tilde{\Theta}_i$. 
	%Then, by \eqref{convex-strong2}, $V(\bm{x})$ further satisfies that 
	%\begin{align}\label{V_bb}
	%\sum_{i=1}^N \frac{\theta_i}{2}\|x_i-x^*\|^2  \leq 	V(\bm{x}) \leq \sum_{i=1}^N \frac{\Theta_i}{2}\|x_i-x^*\|^2.
	%\end{align}
	
	In this case, the function $V_1(t)$ can be expressed as $V_1(t) =\frac{1}{2}\bm{\xi}^T(t)\bm{\xi}(t) $. Then, it can be calculated that
	\begin{align}\label{V1}
		\dot{V}_1 = -\bm{\xi}^T \overline{B}^T \overline{P} \overline{B} \bm{\xi} - \bm{\xi}^T \overline{B}^T (\overline{P} \tilde{y}_x + \overline{Q} \tilde{y}_{\lambda})
	\end{align}  
	with $\tilde{y}_x =[a_iy_{x_i}]_{i\in \mathcal{V} }$, $\tilde{y}_{\lambda} =[a_iy_{\lambda_i}]_{i\in \mathcal{V} }$ and $\overline{Q} = \text{blkdiag}(Q_1,\cdots,Q_N)$. Since 
	$\bm{x}$ is bounded and $y_i(t) = y_i(0)e{-a_it}$, then there exists a constant $\rho>0$ such that $|\bm{\xi}^T \overline{B} (\overline{P} \tilde{y}_x + \overline{Q} \tilde{y}_{\lambda})|\leq \rho e^{-\underline{a} t}$. Let $M = \overline{B}^T \overline{P} \overline{B}$. With Lemma \ref{lem-opt}, we can conclude that $M\neq 0$ and $\mathcal{N}(M) = \mathcal{N}(\overline{B})$, which implies that $\bm{\xi}  \perp \mathcal{N}(M)$. When $\bm{\xi}\neq 0$, it holds that $\bm{\xi}^T \overline{B}^T \overline{P} \overline{B} \bm{\xi} \geq \lambda_{2}(M)\bm{\xi}^T\bm{\xi}$ by Lemma \ref{lem-lam}, where $\lambda_{2}(M)>0$ is a positive continuous function dependent on $\bm{x}$ according to Lemma \ref{lem-lam2}. As $\bm{x}\in U$, then there exist a lower bound $\lambda_0$ such that $\lambda_{2}(M) \geq \lambda_0, \forall \bm{x}\in U$. As a result, \eqref{V1} can be further relaxed as
	\begin{align}\label{V1-relax}
		\dot{V}_1 \leq -2\lambda_0V_1 + \rho e^{-\underline{a} t}. 
	\end{align}  
	By the comparison theorem, $V_1$ satisfies the inequality 
	\begin{align*}
		V_1 &\leq V_1(0) e^{-2\lambda_0 t} + \rho \int_{0}^{t} e^{-2\lambda_0(t-\tau)}e^{-\underline{a} \tau}d \tau   \\
	%	&= V_1(0) \text{exp}(-2\lambda_0 t) + \rho \text{exp}(-2\lambda_0t)(\frac{\text{exp}((2\lambda_0-\underline{a})t)-1}{2\lambda_0-\underline{a}}) \\
		&= (V_1(0)-\frac{\rho}{2\lambda_0-\underline{a}}) e^{-2\lambda_0 t} +  \frac{\rho e^{-\underline{a}t}}{2\lambda_0-\underline{a}} 
	\end{align*}
	for which $\lambda_0\neq \underline{a}$ without loss of generality, which indicates that $V_1$ converges to zero exponentially. Let $\bar{x}= \bm{1}_N^T\bm{x}/N $ and $\bar{z}=(\bar{x},\lambda)$. Then, it can be easily verified $x_i$ converges to $\bar{x}$ exponentially since there exists a constant $\tau>0$ such that $I_N-\bm{1}_N\bm{1}_N^T/N\leq \tau B^TB$. One can suppose that $\|\bm{z}-\bm{1}_N\otimes \bar{z}\|\leq C_1 e^{-\nu_1 t}$ for some $C_1, \nu_1>0$. By applying the mean-value theorem to expand $\nabla_z \mathcal{L}(\bar{z})$ around the optimal primal-dual point $z^*=(x^*,\lambda^*)$, we have
	\begin{align}\label{mean}
	 	\nabla \mathcal{L}(\bar{z}) = (\nabla^2\mathcal{L}(\eta(\bar{z})))^{-1}(\bar{z}-z^*),
	 \end{align}
	where $\eta(\bar{z})$ is a convex combination of $\bar{z}$ and $z^*$. Since $\nabla^2 \mathcal{L}(\theta(\bar{z}))$ is positive definite and $\bar{x}$ is bounded due to $\bm{x}\in U$, then there exists a constant $\tau_1$ such that $\|\nabla^2\mathcal{L}(\eta(\bar{z}))\|^{-1} \leq \tau_1$. Therefore, we can further derive from \eqref{mean} that
	\begin{align}\label{mean2}
	\|\bar{z}-z^*\| = \| (\nabla^2\mathcal{L}(\eta(\bar{z})))^{-1} \nabla \mathcal{L}(\bar{z}) \| \leq \tau_1 \|\nabla_z \mathcal{L}(\bar{z})\|. 
    \end{align}
	Moreover, we have $\nabla_{\lambda} \mathcal{L}(\bar{z}) = [\nabla_{\lambda_i} \mathcal{L}_i(z_i)]_{i \in \mathcal{V}} $ and 
	\begin{align*}
	& \|\nabla_x \mathcal{L}(\bar{z}) - \sum_{i\in \mathcal{V}}\nabla_{x_i} \mathcal{L}(z_i)\| = \| \sum_{i\in \mathcal{V}} (\nabla f_i(\bar{x})- \nabla f_i(x_i))\|\\
	& = \| \sum_{i\in \mathcal{V}} H_i(\eta_i) (\bar{x}- x_i)\| \\
	& \leq \sum_{i\in \mathcal{V}}\|  H_i(\eta_i)\|\| \bar{x}- x_i\| \leq \tau_2C_1 e^{-\nu_1 t}
	\end{align*}
	where the second equality holds by applying the mean-value theorem around $\bar{x}$, and the positive constant $\tau_2$ exists in the last inequality since $\bar{x}$ and $x_i$ are bounded and $f_i$ is twice continuously differentiable. As both $\sum_{i=1}^{N}\nabla_{x_i} \mathcal{L}(z_i)$ and $\nabla_{\lambda_i} \mathcal{L}_i(z_i)$ converge to zero exponentially, so does $\nabla \mathcal{L}(\bar{z})$. Hence, by \eqref{mean2}, $\bar{z}$ converges to $z^*$ exponentially. Then, $\bm{z}(t)$ converges to $\overline{\mathcal{Z}}^*$ exponentially.

	(2) With $\eta=0$ and $\alpha_i\in [0 \ 1), \alpha_{ij}=1$, the subsystem \eqref{EZGS-y} is FT stable with settling time bounded by $T_0$ for all agents, and then $x_i$ follows the following discontinuous dynamics
	\begin{align}\label{ZGS-xsign}
		\dot{x}_i(t) &= -P_i(x_i)\sum_{j\in \mathcal{N}_i}\text{sgn}(x_i-x_j), \forall t \geq T_0. 
	\end{align}
	By the proof of Theorem \ref{thm-ZGS}, $\bm{x}$ belongs to a compact set $\Phi$. Therefore, the right-hand side of \eqref{ZGS-xsign} is bounded and the its Filippov solution always exists. Considering the function $V_1=\|\overline{B}\bm{x}\|_1$, based on the nonsmooth analysis \cite[Prop. 10]{Cortes2008}, there exists $\nu \in \mathcal{F}[\text{sgn}\circ\overline{B} ](\bm{x}) = \partial V_1(\bm{x})$ such that 
	\begin{align}\label{V1-sign}
		\dot{V}_1(t) & \overset{a.e.}{=} -\nu^T \overline{B}^T \overline{P} \overline{B} \nu = -\nu^T M  \nu,
	\end{align}
	where $\mathcal{F}[\psi](\cdot)$ denotes the Filippov set-vauled map of $\psi$ and $\partial V_1(\cdot)$ represents the generalized gradient of $V_1$ \cite{Cortes2008}. As $M$ is positive semi-definite, it holds that
	\begin{align}
		\nu^T M \nu \geq \lambda_2(M)\|\nu - \mathcal{P}_{\mathcal{N}(M)}(\nu)\|_2^2,
	\end{align}
    where $\mathcal{P}_{\mathcal{N}(M)}(\nu)$ denotes the projection of $\nu$ onto $\mathcal{N}(M)$.
	By Lemma \ref{lem-opt},  $\mathcal{N}(M)=\mathcal{N}(\overline{B})$. Moreover, 
	\begin{align*}
		\|\bm{\xi}\|_1 = \nu^T \bm{\xi} &= (\mathcal{P}_{\mathcal{N}(\overline{B})}(\nu)+\nu-\mathcal{P}_{\mathcal{N}(\overline{B})}(\nu) )\bm{\xi}\\
		&= (\nu-\mathcal{P}_{\mathcal{N}(\overline{B})}(\nu) )\bm{\xi}\leq \|\nu-\mathcal{P}_{\mathcal{N}(\overline{B})}(\nu) \|_2\|\bm{\xi}\|_2.
	\end{align*}
	Since $\|\bm{\xi}\|_1 \geq \|\bm{\xi}\|_2$ when $\bm{\xi}\neq 0$. Hence, $\|\nu-\mathcal{P}_{\mathcal{N}(\overline{B})}(\nu) \|_2\geq 1$ when $\bm{\xi}\neq 0$. Therefore, $\dot{V}_1 \leq -\lambda_{2}(M)$ when $\bm{\xi}\neq 0$. By Lemma \ref{lem-lam2}, $\lambda_2(M)$ is a positive continuous function of $\bm{x}$. Then, there must exist a positive scalar $\lambda_0$ such that $\lambda_2(M)\geq \lambda_0, \forall \bm{x} \in \Phi$. As a result, $\dot{V}_1 \leq -\lambda_{2}(M)$ when $V_1\neq 0$. As $V_1$ is positive definite, then $V_1$ will converge to zero in FT $T = T_0+\frac{V_1(T)}{\lambda_0}$, which means that $\bm{x}$ reaches $\mathcal{A}$ in FT. As $\bm{z}(t) \in \mathcal{S}_0$ for $t\geq T$, then $\bm{z}(t)$ will converge to $(\mathcal{A}\times \mathbb{R}^m)\cap\mathcal{S}_0 = \overline{\mathcal{Z}}^*$ in FT $T$.
	
	(3) With $\eta=0$ (resp. 1) and $\alpha_i\in [0 \ 1), \alpha_{ij}\in (0 \ 1)$, the subsystem \eqref{EZGS-y} is FT (resp. FxT) stable, and it is suppose that the settling time is bounded by $T_0$ for all agents. According to Theorem \ref{thm-ZGS}, it holds that $z(t) \in \mathcal{S}_0=\mathcal{M} \cap (\mathcal{C}\times \mathbb{R}^m), \forall t\geq T_0$. Then, it remains to show that $\bm{x}(t)$ will converge to $\mathcal{A}$ in FT after $T_0$. When $t\geq T_0$, $x_i$ follows the dynamics \eqref{ZGS-x}. Taking the derivative of $V_1$ along \eqref{ZGS-x} gives
	\begin{align}\label{V1-FT}
		\dot{V}_1 = -\chi^T(\bm{\xi}) \overline{B}^T \overline{P} \overline{B} \chi(\bm{\xi}) \leq - \lambda_{\chi}(M) \chi^T(\bm{\xi}) \chi(\bm{\xi}),
	\end{align} 
	where the inequality is due to Lemma \ref{lem-lam} with $\lambda_{\chi}(M)>0$ being a positive continuous function of $\bm{x}$ by Lemma \ref{lem-lam2}. According to the proof of Theorem \ref{thm-ZGS}, $\bm{x}$ belongs to a compact set and hence there exists $\lambda_0>0$ (here we use the same symbol as case (1) if no confusion arise) such that $\lambda_{\chi}(M)\geq \lambda_0$. Then, referring to \eqref{V1-FT}, we have 
%	\begin{align}\label{V1-FT2}
%		\dot{V}_1 \leq  - \lambda_0 \chi^T(\bm{\xi}) \chi(\bm{\xi})=- \lambda_0 (\text{sgn}^{[\bar{\alpha}]}(\bm{\xi}))^T (\text{sgn}^{[\bar{\alpha}]}(\bm{\xi})).
%	\end{align} 
    \begin{align}\label{V1-FT2}
    	\dot{V}_1 &\leq  - \lambda_0 \chi^T(\bm{\xi}) \chi(\bm{\xi})=- \lambda_0   \|\text{sgn}^{[\bar{\alpha}]}(\bm{\xi})+\eta\text{sgn}^{[\bar{\beta}]}(\bm{\xi})\|^2. \nonumber
    \end{align} 
	Then, similar to the proof of \cite[Therorem 1]{Shi2020auto}, one can show that $\bm{\xi}$ converges to zero in FT $T>T_0$. Then, following case (2), $\bm{z}(t)$ converges to $\overline{\mathcal{Z}}^*$ in FT $T$. 
	
	(4) With $\eta=1$ and $\alpha_i\in [0 \ 1), \alpha_{ij}\in (0 \ 1)$, if there exists a scalar $\underline{\lambda}_0$ such that $\lambda_{\chi}(M(\bm{x}))\geq \underline{\lambda}_0, \forall \bm{x}\in \mathbb{R}^{nN}$, \eqref{V1-FT2} holds and we have $\lambda_0\geq \underline{\lambda}_0$. Then, similar to the proof of \cite[Therorem 1-2)]{Shi2020auto}, one can show that $\bm{\xi}$ converges to zero in FxT $T>T_0$. Hence, $\bm{z}(t)$ converges to $\overline{\mathcal{Z}}^*$ in FxT $T$. The proof is complete. 
	%The proof is complete.

%	(4) With $\eta=1$ and $\alpha_i\in [0 \ 1), \alpha_{ij}\in (0 \ 1)$, the subsystem \eqref{EZGS-y} is FxT stable with the settling time $T_0$ independent of the initial states. Then, similar to the proof of case (2), it can be shown that there exists a positive scalar $\lambda_0$ such that
%	\begin{align*}
%		\dot{V}_1 &\leq  - \lambda_0 \chi^T(\bm{\xi}) \chi(\bm{\xi}) \\
%		&=- \lambda_0  (\text{sgn}^{[\bar{\alpha}]}(\bm{\xi})+\text{sgn}^{[\bar{\beta}]}(\bm{\xi}))^T (\text{sgn}^{[\bar{\alpha}]}(\bm{\xi}+\text{sgn}^{[\bar{\beta}]}(\bm{\xi}))).
%	\end{align*} 
%	Then, similar to the proof of \cite[Therorem (2)]{Shi2020auto}, one can show that $\bm{\xi}$ converges to zero in FxT $T>T_0$ by applying Lemma \ref{lem-FTC}. Hence, $z(t)$ converges to $\overline{\mathcal{Z}}^*$ in FxT $T$. The proof is complete.

\subsection{Proof of Theorem \ref{thm3}}\label{App-C}
Consider the function $V_i(t) = \frac{1}{2}y_i^Ty_i$ for the subsystem \eqref{EZGS-y}. The derivative of $V_i(t)$ satisfies 
\begin{align}
	\dot{V}_i = -2d V_i  -2 \frac{\dot{\mu}(t;T_0)}{\mu(t;T_0)}V_i, 
\end{align}
implying that $y_i$ converges to 0 in PT $T_0$. Then, by \eqref{PT-Vderi11}, we have 
\begin{align}\label{u1}
	\|y_i(t)\|  \leq  \mu(t;T_0)^{-1}e^{-dt} \|y_i(0)\|, \ \forall t\in [0, T_0). 
\end{align}
Meanwhile, for $t\in [0, T_0)$, with \eqref{u1}, we have 
\begin{align*}
	\|g_i(y_i(t),t)\|  &= \|(c_i + \frac{ \dot{\mu}(t; T_0)}{ \mu(t; T_0)})y_i(t) \| \\
	& = \|(c_i + \frac{ \dot{\mu}(t; T_0)}{ \mu(t; T_0)}) \mu(t;T_0)^{-1}e^{-c_it} \|y_i(0)\|\\
	& \leq (c_i (\frac{T_0-t}{T_0})^h + \frac{h}{T_0}(\frac{T_0-t}{T_0})^{h-1}) e^{-c_it} \|y_i(0) \|\\
	%	& \leq |c_i \mu(t;T_0)^{-1} + \frac{h}{T_0}(\frac{T_0-t}{T_0})^{h-1}|e^{-c_it} \|y_i(0) \|\\
	& \leq (c_i +  \frac{h}{T_0})e^{-c_it} \|y_i(0) \|.
\end{align*}
Hence, $\|g_i(y_i(t),t)\| \in L_{\infty}$ and $\lim_{t\rightarrow T_0^{-}}g_i(y_i(t),t) = 0$, where $L_{\infty}= \{s(t): s:\mathbb{R}_+ \rightarrow \mathbb{R}, \text{sup}_{t\in \mathbb{R}_+}|s(t)|<\infty \}$. Combining with $g_i(y_i(t),t)\equiv0, \forall t\geq T_0$, it indicates that the right-hand side of \eqref{EZGS-y} is continuous and uniformly bounded over $[0, \infty)$. Similar to the proof of Theorem \ref{thm-ZGS}, one  can show that $z(t) \in \mathcal{S}_0=\mathcal{M} \cap (\mathcal{C}\times \mathbb{R}^m), \forall t\geq T_0$.

In the next, we will show $\bm{x}$ is uniformly bounded over $[0, \infty)$. Consider the same function $V(\bm{x})$ expressed by \eqref{V0}, referring to \eqref{V_relax}, for any $t\geq 0$, we have 
\begin{align*}
	\dot{V}(\bm{x}) \leq &  \sum_{i \in \mathcal{V}}  (c_i +  \frac{h}{T_0})e^{-c_it} \|y_{x_i}(0)\| \|e_i(t)\|  \\
	 \leq & \bar{w} e^{-\underline{c} t} \sqrt{  2N/\underline{\theta}}\sqrt{V},
	%		=&-\frac{1}{2}\sum_{i \in \mathcal{V} }\sum_{j\in \mathcal{N}_i}\|x_j-x_i\|^2 - \sum_{i=1}^N a_i(x_i-x^*)^Ty_i(t)|_{x_i} \\
	%	& + a_i(y_i(t)|_{\lambda_i})^T(Q_iy_i(t)|_{x_i} -S_i^{-1}y_i(t)|_{\lambda_i})\\
\end{align*}
where $\bar{w} = \max_{i \in \mathcal{V}}\{(c_i +  \frac{h}{T_0}) \|y_i(0)|_{x_i}\|\}$ and $\underline{c} =\min_{i \in \mathcal{V}}\{c_i \} $. Then, it can be calculated that $V$ satisfies $V(\bm{x}(t)) \leq (\sqrt{V_0}+\underline{c}\bar{w}\sqrt{N/(2\underline{\theta}})) \triangleq \varepsilon_2, \forall t\geq 0 $.
%	for any $t\geq 0$
%	\begin{align}
%		V(\bm{x}(t)) \leq (\sqrt{V_0}+\underline{c}\bar{w}\sqrt{\frac{N}{2\underline{\theta}}}) \triangleq \varepsilon_2 .
%	\end{align} 
Then, It indicates that $\overline{U}= \{\bm{x}\in \mathbb{R}^{Nn}: 	V(\bm{x}) \leq \varepsilon_2\}$ and $\overline{U}_i= \{x_i\in \mathbb{R}^n: f_i(x^*)-f_i(x_i)-\nabla^T f_i(x_i) (x^*-x_i) \leq \varepsilon_2\}$ are compact. We further have $\|\bm{e}\|\leq \sqrt{2\varepsilon_2/\underline{\theta}}$ by \eqref{V_lb0}. Since $f_i$ is twice continuously differentiable, then there exists $\overline{\Theta}_i >0$ such that $H_i(x_i)\leq \overline{\Theta}_i I$ for any $i\in \mathcal{V}$, which indicates that $\theta_i \leq \|H_i\|\leq \overline{\Theta}_i$.

Then, it remains to show that $\bm{x}(t)$ will converge to $\mathcal{A}$ at PT $T \geq T_0$. With the same notation $\bm{\xi}$ as in the proof of Theorem \ref{thm2}, it holds that 
\begin{align}\label{dxi}
	\dot{\bm{\xi}}= - (d + \kappa \frac{ \dot{\mu}(t; T)}{ \mu(t;T)}) B \overline{P} B^T \bm{\xi}- B \overline{P} g_x(y,t), 
\end{align}
where $g_x(y,t) =[g_{x_i}(y_i,t)]_{i \in \mathcal{V}}$. Take the Lyapunov candidate $V_1(t) = \frac{1}{2}\bm{\xi}^T\bm{\xi}$. When $\bm{\xi}\neq 0$, the derivative of $V(t)$ is given by 
\begin{align}\label{dV1}
	\dot{V_1} = - (d + \kappa \frac{ \dot{\mu}(t; T)}{ \mu(t; T)}) \bm{\xi}^T B \overline{P} B^T \bm{\xi}^T- \bm{\xi}^T B \overline{P} g_x(y,t). 
	%\\ &\leq - (d + \kappa \frac{ \dot{\mu}(t; T)}{ \mu(t; T)}) \lambda_2(M) \bm{\xi}^T\bm{\xi}^T + \bar{d} \mu^{-1+\frac{1}{h}}(t; T),  
\end{align}
Denote $\mu(t) = \mu(t; T)$. When $T=T_0$, applying Young’s inequality, for any $\epsilon_1, \epsilon_2>0$, it can be derived that 
\begin{align*}
	\dot{V_1}& =  \frac{ \dot{\mu}}{ \mu} ( -\kappa\bm{\xi}^T B \overline{P} B^T \bm{\xi}^T + \bm{\xi}^T B \overline{P} \bm{y}) \\
	 & \qquad \qquad \qquad\qquad   + d (-\bm{\xi}^T B \overline{P} B^T \bm{\xi}^T +\bm{\xi}^T B \overline{P}\bm{y})\\
	 &\leq  \frac{ \dot{\mu}}{ \mu} ( -(\kappa-\frac{\epsilon_1}{2})\bm{\xi}^T B \overline{P} B^T \bm{\xi}^T + \frac{1}{2\epsilon_1}\bm{y}^T \overline{P} \bm{y}) \\
	 & \qquad \qquad \qquad   + d (-(1- \frac{\epsilon_2}{2})\bm{\xi}^T B \overline{P} B^T \bm{\xi}^T +\frac{1}{2\epsilon_2}\bm{y}^T \overline{P} \bm{y})\\
	 &\leq -(d (1- \frac{\epsilon_2}{2}) +(\kappa-\frac{\epsilon_1}{2})\frac{ \dot{\mu}}{\mu}) \lambda_{2}(M) \bm{\xi}^T\bm{\xi}^T + \bar{d} \mu^{-2+\frac{1}{h}}
	%\\ &\leq - (d + \kappa \frac{ \dot{\mu}(t; T)}{ \mu(t; T)}) \lambda_2(M) \bm{\xi}^T\bm{\xi}^T + \bar{d} \mu^{-1+\frac{1}{h}}(t; T),  
\end{align*}
where the inequality holds due to Lemma \ref{lem-lam}, \eqref{u1} and the fact that there exists a bound $\bar{d}>0$ depending on $\epsilon_1, \epsilon_2$ such that 
\begin{align*}
\frac{\dot{\mu}}{2\epsilon_1\mu}\bm{y}^T \overline{P} \bm{y}+ \frac{d}{2\epsilon_2}\bm{y}^T \overline{P} \bm{y}\leq \bar{d} \mu^{-2+\frac{1}{h}}
\end{align*}
as $\frac{1}{h}<1$ and $\bm{\xi}$ is bounded. Since $\bm{x}$ belongs to a compact set $\overline{U}$, then there exists a positive scalar $\lambda_0$ such that $\lambda_2(M) \geq \lambda_0$ due to Lemma \ref{lem-lam2}. Hence, when $\kappa > 1/\lambda_0, \epsilon_1\in (0, 2(\kappa- 1/\lambda_0)) $ and $\epsilon_2 \in (0, 2)$, we further have 
\begin{align}\label{PT-V2}
	\dot{V_1} \leq - (2\lambda_0 \tilde{d} +2 \frac{ \dot{\mu}}{ \mu}) V_1+  \bar{d} \mu^{-2+\frac{1}{h}}.  
\end{align} 
Then, by applying Lemma \ref{lem-PT}, we have 
\begin{align}\label{PT-Vderi2}
	V_1(t)  \leq \mu^{-2}e^{-2\lambda_0 \tilde{d} t}V_1(0)+\bar{d}T \mu^{-2} \text{ln}(\frac{T}{T-t}). 
\end{align}
As $\lim_{ t\rightarrow T}\mu^{-2}(t) \text{ln}(\frac{T}{T-t})= 0$, we can conclude that $V_1(t)$ $(\text{i.e.}, \bm{\xi})$ converges to 0 in PT $T$ and $V_1(t)=0$ $(\text{i.e.}, \bm{\xi}=0), \forall t\geq T$, i.e., $\bm{x}(t) \in \mathcal{A}, \forall t\geq T$. Together with $\bm{z}(t) \in \mathcal{S}_0=\mathcal{M} \cap (\mathcal{C}\times \mathbb{R}^m), \forall t\geq T_0$, we can conclude that $\bm{z}(t) \in \overline{\mathcal{Z}}^*, \forall t \geq T$.

Then, we will show that the right side of EZGS \eqref{EZGS} is continuous and uniformly bounded over $[0, \infty)$. Since $\theta_i \leq \|H_i\|\leq \overline{\Theta}_i, \forall i\in \mathcal{V}$, there exists a bounded scalar $\varrho$ only determined by $\bm{z}(0)$ such that $\|(\nabla^2 \mathcal{L}_i(z_i))^{-1}\|\leq \varrho, \forall i\in \mathcal{V}, t\geq 0$. Then, by \eqref{EZGS}, we have 
	\begin{align}\small
		\|\dot{z}_i(t)\|& \leq \|(\nabla^2 \mathcal{L}_i(z_i))^{-1}\|(\|g_i(y_i,t)\| + k(t)\|\sum_{j\in 
			\mathcal{N}_i} \chi_{ij}(x_i,x_j) \|) \nonumber \\
		& \leq \varrho\|g_i(y_i,t)\| + \varrho(d + \kappa \frac{ \dot{\mu}}{ \mu}) \|\overline{B}^T\|\|\bm{\xi}(t)\|, \label{dzi}
	%	& \leq \varrho (c_i +  \frac{h}{T})\|y_i(0) \| + \varrho(d + \frac{\kappa h}{T-T_0}) \|L\| \sqrt{  \frac{2\varepsilon_2}{\underline{\theta}}}, 
	\end{align}
    By \eqref{PT-Vderi2}, it holds that 
    \begin{align}
    \lim_{ t\rightarrow T^-}	\frac{ \dot{\mu}}{ \mu} \|\bm{\xi}(t)\| = \lim_{ t\rightarrow T^-}\frac{h}{T} \mu^{\frac{1}{h}} \sqrt{2V_i(t)}=0. 
    \end{align}
Combining with \eqref{u1}, one can conclude that $ \lim_{ t\rightarrow T}	\|\dot{\bm{z}}_i(t)\| = 0, \forall i\in \mathcal{V}$, equal to the control law of \eqref{EZGS-z} at $T$, indicating that the control law of \eqref{EZGS} is continuous at $t=T$. Moreover, the control law of \eqref{EZGS} is uniformly indicating by \eqref{u1} and \eqref{dzi}. 

When $T>T_0$, it is obvious that $u_i(t)$ is bounded over $[0, T_0]$ as $\bm{x}(t)$ is bounded due to $\frac{\dot{\mu}(t; T)}{ \mu(t; T)} \in [\frac{h}{T},\frac{h}{T-T_0}] $. For $t\in [T_0, T)$, $g_i(y,t)=0$. By \eqref{dV1} with $\kappa \geq \frac{1}{\lambda_0}$, we get that 
\begin{align*}
	\dot{V_1} = - (d + \kappa \frac{ \dot{\mu}(t; T)}{ \mu(t; T)}) \bm{\xi}^T B \overline{P} B^T \bm{\xi}^T  \leq - (2\lambda_0 d +2 \frac{ \dot{\mu}}{ \mu}) V_1. 
\end{align*}
By Lemma \ref{lem-PT}, it yields that 
\begin{align*}
	V_1(t)  \leq \mu^{-2}e^{-2\lambda_0 d t}V_1(0), 
\end{align*}
implying that 
\begin{align*}
		\|\bm{\xi}(t)\|  \leq  \mu(t-T_0;T-T_0)^{-1}e^{-\lambda_0 (t-T_0)} \|\bm{\xi}(T_0)\|, \forall t\in [T_0, T)
\end{align*}
	since $\frac{ \dot{\mu}(t; T)}{ \mu(t; T)} = \frac{ \dot{\mu}(t-T_0;T-T_0)}{ \mu(t-T_0;T-T_0)}$. Then, similar to \eqref{dzi}, one can derive that 
	\begin{align*}\tiny
		\|\dot{\bm{z}}(t)\|
		& \leq \varrho (d + \kappa \frac{ \dot{\mu}(t; T)}{ \mu(t; T)}) \|\overline{B}^T\|\|\bm{\xi}(t)\|\\
		& \leq \varrho|d + \frac{\kappa h}{T-T_0}(\frac{T-t}{T-T_0})^{h-1}| \|\overline{B}\|\|\bm{\xi}(T_0)\|\\
		& \leq \varrho(d +  \frac{\kappa h}{T-T_0}) \|\overline{B}\| \|\bm{\xi}(T_0)\|.
	\end{align*}
	Hence, $\|\dot{\bm{z}}(t)\| \in L_{\infty}$ over $[T_0, T)$. In addition, as $\lim_{t\rightarrow T^{-}}\|\dot{\bm{z}}(t)\|=0$ equal to the control law at $T$, it indicates that the control law of \eqref{EZGS} is continuous at $t=T$. Moreover, $\dot{\bm{z}} = 0, \forall t\geq T$ by the definition of $\mu(t;T)$. Then, one can conclude that the control law of EZGS \eqref{EZGS} is continuous and uniformly bounded over $[0, \infty)$.

\bibliographystyle{ieeetr}
\end{document}